\documentclass{elsart}

\usepackage{amssymb,amsmath}

\newenvironment{proof}{\begin{pf}}{\qed\end{pf}}

\newtheorem{algo}{Algorithm}

\newcommand\bN{\mathbb N}
\newcommand\bQ{\mathbb Q}
\newcommand\bR{\mathbb R}

\newcommand{\im}{\operatorname{im}}
\newcommand{\redu}{\operatorname{red}}
\newcommand{\field}{K} 
\newcommand{\ann}{\operatorname{ann}}
\newcommand{\cF}{{\mathcal F}}
\newcommand{\cG}{{\mathcal G}}
\newcommand{\cH}{{\mathcal H}}
\newcommand{\cI}{{\mathcal I}}
\newcommand{\cL}{{\mathcal L}}
\newcommand{\cU}{{\mathcal U}}
\newcommand{\cV}{{\mathcal V}}
\newcommand{\cX}{{\mathcal X}}
\newcommand{\dlAnnF}{I_F^\sadj\big|_{\partial_t=\partial_\ell}}
\newcommand{\drAnnG}{I_G\big|_{\partial_t=\partial_r}}
\newcommand{\dlAnnFt}{\cI_F^\sadj\big|_{\partial_t=\partial_\ell}}
\newcommand{\drAnnGt}{\cI_G\big|_{\partial_t=\partial_r}}
\newcommand{\FplusG}{I_F^\sadj\otimes_{W_p[t]}W_{p,t}^{\vphantom{\sadj}}
+W_{p,t}^\sadj\otimes_{W_p[t]}I_G}
\newcommand{\FplusGt}{\cI_F^\sadj\otimes_{W_p(t)}W_{p,t}(t)^{\vphantom{\sadj}}
+W_{p,t}(t)^\sadj\otimes_{W_p(t)}\cI_G}
\newcommand{\WoW}{W_{p,t}^\sadj\otimes_{W_p[t]}W_{p,t}^{\vphantom{\sadj}}}
\newcommand{\annFxG}{\ann_{\WoW}(F^\sadj\otimes G)}
\newcommand{\annFoG}{\ann_{W_t}(F^\sadj\otimes G)}
\newcommand{\annFoGprime}{\ann_{\Wt}(F^\sadj\otimes G)}
\newcommand{\annFoGt}{\ann_{W_t(t)}(F^\sadj\otimes G)}
\newcommand{\annFoGtprime}{\ann_{\Wt(t)}(F^\sadj\otimes G)}
\newcommand{\Wt}{W'_t}
\newcommand{\rsc}[2]{\left\langle#1,#2\right\rangle} 
\newcommand{\rsp}[2]{\left\langle#1|#2\right\rangle} 

\def\uadj{\star}    
\def\sadj{\diamond} 

\begin{document} 
\begin{frontmatter}
\title{Effective Scalar Products of D-finite Symmetric Functions}
\author{Fr\'ed\'eric Chyzak}
\address{Projet Algorithmes, INRIA Rocquencourt\\{\tt frederic.chyzak@inria.fr}}
\author{Marni Mishna}
\address{LaCIM, Universit\'e du Qu\'ebec \`a Montr\'eal\\{\tt mishna@math.uqam.ca}}
\author{Bruno Salvy}
\address{Projet Algorithmes, INRIA Rocquencourt\\{\tt bruno.salvy@inria.fr}}

\begin{abstract}
Many combinatorial generating functions can be expressed as
combinations of symmetric functions, or extracted as sub-series and
specializations from such combinations.  Gessel has outlined a large
class of symmetric functions for which the resulting generating
functions are D-finite.  We extend Gessel's work by providing
algorithms that compute differential equations these generating
functions satisfy in the case they are given as a scalar product of
symmetric functions in Gessel's class.  Examples of applications to
$k$-regular graphs and Young tableaux with repeated entries are given.
Asymptotic estimates are a natural application of our method, which we
illustrate on the same model of Young tableaux.  We also derive a
seemingly new formula for the Kronecker product of the sum of Schur
functions with itself.  (This article completes the extended abstract
published in the proceedings of FPSAC'02 under the title ``Effective
D-Finite Symmetric Functions''.)
\end{abstract}
\end{frontmatter}

\section*{Introduction} 

A power series in one variable is called differentiably finite, or
simply D-finite, when it is solution of a linear differential equation
with polynomial coefficients. This differential equation turns out to
be a convenient data structure for extracting information related to
the series and many algorithms operate directly on this differential
equation. In particular, the class of univariate D-finite power series
is closed under sum, product, Hadamard product, and Borel transform,
among other operations, and algorithms computing the corresponding
differential equations are known (see for
instance~\cite{Stanley99}). Moreover, the coefficient sequence of a
univariate D-finite power series satisfies a linear recurrence, which
makes it possible to compute many terms of the sequence efficiently.
These closure properties are implemented in computer algebra
systems~\cite{Mallinger96,SaZi94}. Also, the mere knowledge that a
series is D-finite gives information concerning its asymptotic
behavior. Thus, whether it be for algorithmic or theoretical reasons,
it is often important to know whether a given series is D-finite or
not, and it is useful to compute the corresponding differential
equation when possible.

D-finiteness extends to power series in several variables: a power
series is called D-finite when the vector space spanned by the series
and its derivatives is finite-dimensional. Again, this class enjoys
many closure properties and algorithms are available for computing the
systems of linear differential equations generating the corresponding
operator ideals~\cite{Chyzak98,ChSa98}. Algorithmically, the key tool
is provided by Gr\"obner bases in rings of linear differential
operators and an implementation is available in Chyzak's {\tt Mgfun}
package\footnote{This package is part of the {\tt algolib} library
available at {\tt
http://\discretionary{}{}{}algo\discretionary{}{}{}.inria\discretionary{}{}{}.fr/\discretionary{}{}{}packages/}.}. An
additional, very important closure operation on multivariate D-finite
power series is definite integration. It can be computed by an
algorithm called {\em creative telescoping}, due to
Zeilberger~\cite{Zeilberger91b}. Again, this method takes as input
(linear) differential operators and outputs differential operators (in
fewer variables) satisfied by the definite integral. It turns out that
the algorithmic realization of creative telescoping has several common
features with the algorithms we introduce here.

Beyond the multivariate case, Gessel considered the case of infinitely
many variables and laid the foundations of a theory of D-finiteness
for symmetric functions~\cite{Gessel90}.  He defines a notion of
D-finite symmetric series and obtains several closure properties. The
motivation for studying D-finite symmetric series is that new closure
properties occur and can be exploited to derive the D-finiteness of
usual multivariate or univariate power series. Thus, the main
application of~\cite{Gessel90} is a proof of the D-finiteness for
several combinatorial counting functions. This is achieved by
describing the counting functions as combinations of coefficients of
D-finite symmetric series, which can then be computed by way of a
scalar product of symmetric functions. Under certain conditions, the
scalar product of symmetric functions depending on extra parameters is
D-finite in those parameters, where D-finiteness is that of (usual)
multivariate power series. Most of Gessel's proofs are not
constructive. In this article, we give algorithms that compute the
resulting systems of differential equations for the scalar product
operation. Besides Gessel's work, these algorithms are inspired by
methods used by Goulden, Jackson, and Reilly
in~\cite{GoJa83,GoJaRe83}. Finally, Gr\"obner bases are used to help
make these methods into algorithms. One outcome is a simplification of
the original techniques of~\cite{GoJa83,GoJaRe83}.

Considering some enumerative combinatorial problem of a symmetric
flavor and parameterized by a discrete parameter (denoted by~$k$ in
the examples below), it is often so that the enumeration is solved by
first forming a scalar product of two symmetric functions in
$k$~variables.  Moreover, in the examples envisioned (the enumeration
of $k$-regular graphs, of $k$-uniform tableaux, etc.), this scalar
product is the specialization to $k$~variables of a scalar product
between two ``closed form'' symmetric functions in infinitely many
variables.  Both symmetric functions are sufficiently well-behaved
that nice ``closed forms'' are obtained under specialization, leading
to descriptions in terms of linear differential operators that are
easy to derive.  This nice behavior is well exemplified by
Eq.~\eqref{eq:kregp} and Eq.~\eqref{eq:kunif} below and is what
delimits the scope of our method in applications.

Additionally, our method extends to other scalar products whose
associated adjunctions satisfy a certain condition of preservation of
degree (see Section~\ref{sec:other-scalar-products}), as well as to
the Kronecker product of symmetric functions (see
Section~\ref{sec:kronecker-product}).

A very basic example of application of
our method is the enumeration of labeled graphs.  A finite graph on
$n$~vertices labeled with non-negative integers $i_1,\dots,i_n$, of
respective valencies $v_1,\dots,v_n$, is given as a weight the
monomial~$x_{i_1}^{v_1}\dots x_{i_n}^{v_n}$.  This encoding leads to
generating functions that are symmetric series: the set of all finite
simple graphs is enumerated by the product
\[
G(x)=\sum_{G\in\cG}\prod_{(i,j)\in E(G)}x_ix_j=\prod_{i<j}(1+x_ix_j),
\]
as each edge $(i,j)\in E(G)$ is either in the graph or not.  This
series is obviously invariant under renamings of the~$x_i$'s, which
motivates the involvement of symmetric function theory in the
application.  Finite simple graphs whose vertices all have valency two
are called \emph{2-regular graphs}.  Such a graph contributes to~$G$
by a term of the form $x_{i_1}^2\dotsm x_{i_n}^2$.  Therefore,
extracting the sub-series of~$G$ with same monomials as in the series
expansion of~$\prod_{i\in\bN\setminus\{0\}}(1+x_i^2)$, another
symmetric series, results in the generating series of 2-regular graphs
according to the same encoding.  By symmetry, monomials based on
different sets of indices $i_1,\dots,i_n$ of cardinality~$n$ share the
same coefficient in this extracted series.  In this spirit, it will be
shown in Section~\ref{sec:example} that the number of 2-regular graphs
on $n$~vertices is given as the coefficient of~$t^n$ in the series
\[
G_2(t)=\rsc{\exp\bigl((p_1^2-p_2)/2-p_2^2/4\bigr)}{\exp\bigl(t(p_1^2+p_2)/2\bigr)}.
\]
Here, the scalar product is a scalar product for symmetric functions,
to be defined in the next section; it implements the coefficient
extraction.  The variables $t$, $p_1$, and~$p_2$ can be viewed as
standard variables, although $p_1$ and~$p_2$ will be assigned the
symmetric function interpretation $p_1=x_1+x_2+\dotsb$,
$p_2=x_1^2+x_2^2+\dotsb$.  Our purpose in the present paper is to
describe scalar products of symmetric functions like~$G_2(t)$ by a linear
differential equation.  By our method, Algorithm~1 below calculates
that $G_2(t)$ satisfies the differential equation
\[2(1-t)G'_2(t)-t^2G_2(t)=0,\]
which is easily solved to recover the classical series $G_2(t)=e^{-\frac14t(t+2)}/\sqrt{1-t}$.
More details on this calculation as well as similar examples will be
given in Section~\ref{sec:example}. In general, the derived differential equation will not admit of such a closed form solution. However it is possible to extract asymptotic information on the sequence being enumerated directly from this differential equation. This will be exemplified in Section~\ref{sec:asympt}.

This article is organized as follows. After recalling the necessary
part of Gessel's work in Section~\ref{sec:intro}, we start by focusing
on the special situation when a single argument of the scalar product
depends on extra parameters.  We present an algorithm for computing
the differential equations satisfied by the scalar product in this
case in Section~\ref{sec:alg}. The application to the example of
$k$-regular graphs is detailed in Section~\ref{sec:example}. Then a
special case where the algorithm can be further refined is described in
Section~\ref{sec:hammond}.  We treat a variant of Young tableaux where
each element is repeated $k$~times in Section~\ref{sec:young}. (These
are in bijection with a generalization of involutions~\cite{Knuth70}.)
The general form of the main algorithm, when both arguments depend on
extra parameters, is given in Section~\ref{sec:general}.  Termination
and correctness of the main algorithms are proved in
Section~\ref{sec:proofs}.  Next, in Section~\ref{sec:asympt} we employ
our algorithms to derive asymptotic estimates of the enumerating
sequences of $k$-regular graphs for $k=1,2,3,4$.  Following this
approach of experimental mathematics, we state a conjecture for
general~$k$.  A discussion on several extensions and applications of
the method closes the paper in Section~\ref{sec:concl}, including the
calculation of a seemingly new formula for the Kronecker product of
the sum of all Schur functions with itself.

\section{Symmetric D-finite Functions}\label{sec:intro}

In this section, we recall the facts we need about symmetric functions,
D-finite functions, and symmetric D-finite functions.

\subsection{Symmetric functions}

We first collect basic definitions, notation,
and results of the theory of symmetric functions.  We
refer to~\cite{Macdonald95,Stanley99} for further results.

Symmetric functions are series in the infinite set of variables
$x_1,x_2,\dotsc$ over a field~$\field$ of characteristic~0, subject to
a certain invariance under renumberings of the variables.  The
$\field$-algebra~$\Lambda$ of symmetric functions is formally defined
as follows.  For each positive integer~$m$, the $\field$-vector space
consisting of the polynomials of~$\field[x_1,\dots,x_m]$ that are
fixed under any permutation of the variables is a graded
$\field$-algebra~$G_m$, the algebra of symmetric polynomials in
$m$~variables.  Here the grading is with respect to the total degree
in the $m$~variables and it induces a chain of graded surjective
homomorphisms~$\pi_m$ from~$G_{m+1}$ onto~$G_m$ defined by
setting~$x_{m+1}$ to~0.  Taking the inverse limit (a.k.a.\ projective
limit) of the system $(\{G_m\},\{\pi_m\})$ results in the graded
$\field$-algebra~$\Lambda$ of symmetric functions.  By restriction of
the algebras~$G_m$ and the maps~$\pi_m$ to homogeneous polynomials in
a fixed degree~$n$, the inductive limit becomes a vector
subspace~$\Lambda_n$ of~$\Lambda$.  We have the relation
$\Lambda=\bigoplus_{n\geq0}\Lambda_n$.

We now recall the definitions of the most frequently used bases of the
ring~$\Lambda$ and vector spaces~$\Lambda_n$.  Denote by
$\lambda=(\lambda_1,\dots,\lambda_k)$ a partition of the integer~$n$.
This means that $n=\lambda_1+\dots+\lambda_k$
and~$\lambda_1\ge\dots\ge\lambda_k>0$, which we also
denote~$\lambda\vdash n$.  Alternatively, the power notation
$\lambda=1^{r_1}\dotsm k^{r_k}$ for partitions indicates that
$i$~occurs $r_i$~times in~$\lambda$, for~$i=1,2,\dots,k$.  Partitions
serve as indices for the five principal symmetric function families
that we use:
\begin{itemize}
\item the homogeneous symmetric functions
$h_\lambda=h_{\lambda_1}\dotsm h_{\lambda_k}$, for $h_n$~defined as
the sum of all monomials of degree~$n$ in~$x_1,x_2,\dotsc$, with
possible repetition (i.e., with any non-negative exponents),
\item the elementary symmetric functions
$e_\lambda=e_{\lambda_1}\dotsm e_{\lambda_k}$, for $e_n$~defined as
the sum of all monomials of degree~$n$ in~$x_1,x_2,\dotsc$, with no
possible repetition (i.e., with exponents 0 or~1, exclusively),
\item the power symmetric functions $p_\lambda=p_{\lambda_1}\dotsm
p_{\lambda_k}$, for $p_n$~defined as the sum of the $n$th~power of all
variables,
\item the monomial symmetric functions
$m_\lambda=\sum_\sigma(r_1!\,r_2!\,\dots)^{-1}
x_{\sigma(1)}^{\lambda_1}\dots x_{\sigma(k)}^{\lambda_k}$, where
$\sigma$~ranges over all permutations of the non-negative integers,
\item the Schur symmetric functions $s_\lambda$, whose intuitive
definition is in terms of the representations of the permutation
group~$S_n$, and that can alternatively be defined as the limit
symmetric function when $n$~tends to infinity of the determinant of
the $n\times n$-matrix with $(i,j)$-entry $h_{\lambda_i-i+j}$.
\end{itemize}
When the indices are restricted to all partitions of the same positive
integer~$n$, any of the five families forms a basis for the vector
space of symmetric polynomials of degree~$n$ in~$x_1,x_2,\dotsc$.  On
the other hand, any of the three families indexed by the
integers~$i\in\bN$, $(p_i)$, $(h_i)$, and~$(e_i)$, is algebraically
independent over~$\bQ$ and generates the algebra~$\Lambda$ of
symmetric functions over~$\field$:
$\Lambda=\field[p_1,p_2,\dotsc]=\field[h_1,h_2,\dotsc]=\field[e_1,e_2,\dotsc]$.
In this work, we shall focus on the basis~$(p_i)$, as we shall
endow~$\Lambda$ with a differential structure will regard to the
variables~$p_i$.

Generating series of symmetric functions live in the larger ring of
symmetric series, $\field[t][[p_1,p_2,\dotsc]]$. There, we have the
generating series of homogeneous and elementary functions:
\[
H(t)=\sum_n h_n t^n =\exp\left(\sum_i p_i\frac{t^i}i\right),\qquad
E(t)=\sum_n e_n t^n =\exp\left(\sum_i (-1)^ip_i\frac{t^i}i\right).
\]

\subsection{Scalar product and coefficient extraction}

The ring of symmetric series is endowed with a scalar product defined
as a bilinear symmetric form such that the bases~$(h_\lambda)$
and~$(m_\lambda)$ are dual to each other:
\begin{equation}\label{eq:scalhm}
\rsc{m_\lambda}{h_\mu}=\delta_{\lambda,\mu},
\end{equation}
where $\delta_{\lambda,\mu}$~is~1 if $\lambda=\mu$ and 0 otherwise.

For a partition in power notation, $\lambda=1^{n_1}\dotsm k^{n_k}$,
the normalization constant
\[
z_\lambda:= 1^{n_1}n_1!\dotsm k^{n_k}n_k!
\]
plays the role of the square of a norm of~$p_\lambda$ in the following important
formula:
\begin{equation}\label{eq:scalp}
\rsc{p_\lambda}{p_\mu}=\delta_{\lambda,\mu}z_\lambda.
\end{equation}

The scalar product is a basic tool for coefficient extraction. Indeed,
if we write $F(x_1,x_2,\dotsc)$ in the form~$\sum_\lambda f_\lambda
m_\lambda$, then the coefficient of $x_1^{\lambda_1}\dotsm
x_k^{\lambda_k}$ in~$F$ is $f_\lambda=\rsc F{h_\lambda}$,
by~\eqref{eq:scalhm}.  Moreover, when $\lambda=1^n$, the identity
$h_{1^n}=p_{1^n}$ yields a simple way to compute this coefficient when
$F$ is written in the basis of the~$p$'s:
\begin{thm}[Gessel; Goulden \& Jackson]\label{thm:theta}
Let $\theta$ be the $\field$-algebra homomorphism from the algebra of
symmetric functions over~$\field$ to the algebra~$\field[[t]]$ of
formal power series in~$t$ defined by $\theta(p_1)=t$, $\theta(p_n)=0$
for $n>1$. Then if $F$ is a symmetric function,
\[\theta(F)=\sum_{n=0}^\infty a_n\frac{t^n}{n!},\]
where $a_n$ is the coefficient of $x_1\dotsm x_n$ in~$F$. 
\end{thm}

Gessel also provides an analogue for this theorem when
$\lambda=1^n2^m$ and $\lambda=1^n3^m$ \cite[Theorems~2--4]{Gessel90}.
Combinations of other degree patterns quickly become arduous to write
explicitly.

\subsection{Plethysm}

Plethysm is a way to compose symmetric functions, which in the
simplest case, amounts to simply scaling the indices on the power
sums.  This inner law of~$\Lambda$, denoted $u[v]$ for $u,v$
in~$\Lambda$, is, for $w=\sum_\lambda c_\lambda p_\lambda$, defined by
the rules \cite{Stanley99}
\begin{multline*}
p_n[w]=\sum_\lambda c_\lambda p_{n\times\lambda_1}p_{n\times\lambda_2}\ldots,
\\
(\alpha u+\beta v)[w]=\alpha u[w]+\beta v[w],
\qquad
(uv)[w]=u[w]v[w],
\end{multline*}
where $\alpha,\beta$ in~$\field$.  For example, consider that
$w[p_n]=p_n[w]$, and in particular that $p_n[p_m]=p_{n\times m}$.
Thus, we see that when we write $w\in\Lambda$ in the power sum basis
as $w=w(p_1,p_2,\dots,p_k,\ldots)$, the scaling effect appears on the
indices as
\[w[p_n]=w(p_{1\times n},p_{2\times n},\dots,p_{k\times n},\ldots).\]

\subsection{D-finiteness of multivariate series}

Recall that a series $F\in\field[[x_1,\dots,x_n]]$ is {\em
D-finite\/} in $x_1,\dots,x_n$ when the set of all partial derivatives
and their iterates, $\partial^{i_1+\dots+i_n}F/\partial x_1^{i_1}\dotsm
\partial x_n^{i_n}$, spans a finite-dimensional vector space over the
field $\field(x_1,\dots, x_n)$.  A {\em D-finite description\/} of a
series~$F$ is a set of differential equations whose solutions in any
$\field(x_1,\dots, x_n)$-vector space share this property.  A typical
example of such a set is a system of $n$~differential equations of the
form
\[
q_1(x)f(x)+q_2(x)\frac{\partial f}{\partial x_i}(x)+\dots+
q_k(x)\frac{\partial^kf}{\partial x_i^k}(x)=0,
\]
where $i$~ranges over $1,\dots,n$, each~$q_j$ is in $K(x_1,\dots,x_n)$
for $1\leq j\leq k$, and $k$ and~$q_j$ depend on~$i$.  Observe that by
a theorem of Stafford \cite[Chapter~5]{Borel87}, any D-finite
series~$F$ admits a D-finite description consisting of only two
differential equations.  However, we do not know how to benefit from
this theoretical result in our computational setting, and it will be
more efficient to compute in a systematic way with non-minimal sets.

The properties of D-finite series we need here are summarized in the
following theorem.
\begin{thm}\label{thm:dprops}
\begin{enumerate}
\item The set of D-finite power series forms a $\field$-subalgebra of
$\field[[x_1,\dots, x_n]]$ for the usual product of series;
\item If $F$ is D-finite in $x_1,\dots,x_n$ then for any subset
of variables~$x_{i_1},\dots,x_{i_k}$ the specialization of~$F$ at
$x_{i_1}=\dots=x_{i_k}=0$ is D-finite in the remaining variables;
\item If $P$ is a polynomial in $x_1,\dots,x_n$, then\/
$\exp P(x)$ is D-finite in $x_1,\dots,x_n$;  
\item If $F$ and $G$ are D-finite in the variables $x_1,\dots,
x_{m+n}$,
then the Hadamard product $F\odot G$ with respect to the variables
$x_1,\dots,x_n$ is D-finite in $x_1,\dots,x_{m+n}.$
\end{enumerate}
\end{thm}
(Recall that the Hadamard product of two series
$\sum_{\alpha\in\bN^k}a_\alpha u^\alpha
\odot\sum_{\beta\in\bN^k}b_\beta u^\beta$ is $\sum_{\alpha\in\bN^k}a_\alpha
b_\alpha u^\alpha$, where $u^\alpha=u_1^{\alpha_1}\dotsm
u_k^{\alpha_k}$.)

These properties are classical~\cite{Stanley99}. The first three are
elementary, the last one relies on more delicate properties  of
dimension and is due to Lipshitz~\cite{Lipshitz88}.

We note at this point that it is usually simple in applications to
provide a D-finite description for a D-finite function, as the latter
is most often given as a polynomial expression in ``atomic'' D-finite
functions, usually well-known special functions.  Given a table of
atomic D-finite descriptions, one bases on the closure properties of
Theorem~\ref{thm:dprops} above and uses algorithms described
in~\cite{ChSa98} in order to derive a D-finite description for the
whole expression.  In our examples, doing this will be
straightforward since our functions will be exponentials of
polynomials.

\subsection{D-finite symmetric functions}
The definition of D-finiteness for series in an infinite number of
variables is achieved by generalizing Property~(2) in
Theorem~\ref{thm:dprops}: $F\in\field[[x_1,x_2,\dotsc]]$ is called
{\em D-finite\/} in the infinitely many variables~$x_i$ if, for any
choice of a finite set~$S$ of positive integers, the specialization
to~0 of each~$x_i$ for $i$ not in~$S$ results in a power series that
is D-finite, in the classical sense, in the variables~$x_i$ for $i$
in~$S$.  In this case, all the properties in Theorem~\ref{thm:dprops}
hold in the infinite multivariate case.

The definition is then tailored to symmetric series by considering the
algebra of symmetric series as generated over~$\field$ by the set
$\{p_1,p_2,\dotsc\}$: a symmetric series is called {\em D-finite\/}
when it is D-finite in the $p_i$'s.

Property~(4) in Theorem~\ref{thm:dprops} has the following very
important consequence:
\begin{thm}[Gessel]\label{thm:rsc_pres_df}
Let $f$ and $g$ be elements of
$\field[[t_1,\dots,t_k]][[p_1,p_2,\dotsc]]$, D-finite in the $p_i$'s
and $t_j$'s, and suppose that $g$ involves only finitely many of the
$p_i$'s. Then $\rsc fg$ is D-finite in the $t_j$'s provided it is
well-defined as a power series.
\end{thm}

We return to the example of regular graphs given in the
introduction.  We shall see in Section~\ref{sec:example} that the
exponential generating series~$G_2$ of 2-regular graphs is given as an
extraction of coefficients from the generating series~$G$ of all
finite simple graphs in the form~$G_2=\rsc G{\exp(h_2t)}$ and we shall
provide the explicit representations
\[
G=\exp\left(\sum_i(-1)^i\frac{p_i^2-p_{2i}}{2i}\right)
\qquad\text{and}\qquad
h_2=\frac{p_1^2+p_2}2.
\]
Both $G$ and~$\exp(h_2t)$ are clearly D-finite symmetric by the
definition above.  Now, $G_2$~is equal to the scalar product
\[\rsc{\exp\left(\sum_i(-1)^i(p_i^2-p_{2i})/2i\right)}{\exp\bigl(t(p_1^2+p_2)/2\bigr)},
\]
and thus by Theorem~\ref{thm:rsc_pres_df} the resulting power series
is D-finite in~$t$.  Note the effect of the requirement that $g$~be
dependent on finitely many~$p_i$'s only in the theorem---here $\exp
h_2t$~depends on $p_1$ and~$p_2$ only.  As a consequence, the scalar
product extracts those terms from~$G$ that are supported by monomials
in $t$, $p_1$, and~$p_2$ only.  In other words, we can set all~$p_i$'s
to~0 in~$G$ when~$i>2$, which yields
\[
G_2(t)=\rsc{\exp\bigl((p_1^2-p_2)/2-p_2^2/4\bigr)}{\exp\bigl(t(p_1^2+p_2)/2\bigr)}.
\]
This scalar product is between symmetric functions in finitely
many~$p_i$'s.

\subsection{Effective D-finite symmetric closures}
Our work consists in making Theorem~\ref{thm:rsc_pres_df} effective by
giving algorithms for producing linear differential equations
annihilating~$\rsc fg$.  The input to our algorithms consists of
closed forms for~$g$ and the specialization of~$f$ in the finite
number of~$p_i$'s appearing in~$g$, from which generators of ideals of
differential operators which annihilate them can then be computed.

Providing algorithms to manipulate linear differential equations
amounts to making the closure properties of univariate D-finite series
effective; similarly, algorithms operating on systems of linear
differential operators make the closure properties of multivariate
D-finite series effective.  Our title is thus motivated by the fact
that our algorithm makes it possible to compute all the information on
a scalar product that can be predicted from D-finiteness.  Note that
we do not check that the resulting power series is well-defined: our
algorithm merely computes equations that the scalar product series
must satisfy if it is well-defined.

In our examples, we make use of symmetric series that are built by
plethysm. Closure properties are given by Gessel, but our applications
require only a simple consequence of Property~(3) in
Theorem~\ref{thm:dprops}, namely that if $g$~is a polynomial in the
$p_i$'s, then $h[g]$ and $e[g]$ are D-finite for $h=H(1)$
and~$e=E(1)$.

\section{Algorithm for Scalar Product: the Simple Case}\label{sec:alg} 

We proceed to give a new algorithm to compute the differential
equation satisfied by a scalar product of two D-finite symmetric
series under the hypotheses of Theorem~\ref{thm:rsc_pres_df} and with
the additional simplifying condition that only one of the symmetric
series depends on~$t$.  When the number of $t$~variables is~1, the
output is a single differential equation for which existing computer
algebra algorithms might find a closed-form solution. In most cases
however, no such solution exists and we are content with a
differential equation from which useful information can be extracted.

The basic tool we use here is non-commutative Gr\"obner bases in
extensions of Weyl algebras. An introduction to this topic can be
found in~\cite{SaStTa00}.  By~$W_t$, we denote the Weyl algebra
\begin{multline*}
W_t=\field\bigl\langle t_1,\dots,t_k,\partial_{t_1},\dots,\partial_{t_k};\\
[\partial_{t_i},t_j]=\delta_{i,j},
[t_i,t_j]=[\partial_{t_i},\partial_{t_j}]=0,
\ 1\leq i,j\leq k\bigr\rangle,
\end{multline*}
where the bracket~$[a,b]$ denotes $ab-ba$ and $\delta_{i,j}$~is the
Kronecker notation.  This algebra can be identified with the algebra
of linear differential operators with coefficients that are polynomial
in~$t=t_1,\dots,t_k$.  We correspondingly denote~$W_p$ for variables
$p=p_1,\dots,p_n$, as well as $\partial_t$
for~$\partial_{t_1},\dots,\partial_{t_k}$, $\partial_p$
for~$\partial_{p_1},\dots,\partial_{p_n}$, etc.  For the algorithm, we
work in the extension
\[W_{p,t}(t)=\field(t)\otimes_{\field[t]}W_{p,t}\]
of the Weyl algebra~$W_{p,t}=W_p\otimes_\field W_t$ in which the
coefficients of the differential operators are still polynomial in~$p$
but rational in~$t$.  Suppose $F$ and~$G$ belong to $\field[t][[p]]$
and are D-finite symmetric series as in
Theorem~\ref{thm:rsc_pres_df}. In particular, they both satisfy
systems of linear differential equations with polynomial coefficients
from $\field(t)[p]$. We can write these equations as elements of
$W_{p,t}(t)$ acting on~$F$ and~$G$.  The set
$\cI_F=\ann_{W_{p,t}(t)}F$ (resp.~$\cI_G$) of all operators
of~$W_{p,t}(t)$ annihilating~$F$ (resp.~$G$) is then a {\em left\/}
ideal of~$W_{p,t}(t)$. Given as input Gr\"obner bases of $\cI_F$
and~$\cI_G$, our algorithm outputs non-zero elements of the
annihilating left ideal $\ann_{W_t(t)}\rsc FG$.

To combine elements of $\cI_F$ and~$\cI_G$ in a meaningful way we use
the adjunction map, denoted~$\sadj$ here\footnote{Macdonald denotes
the adjunction operator by~$\perp$.}, defined for an operator $P\in
W_p$ by imposing the relation $\rsc{P\cdot F}G=\rsc F{P^\sadj\cdot G}$
for all series $F$ and~$G$.  As a consequence, we have the relation
$(PQ)^\sadj=Q^\sadj P^\sadj$ and the adjoint~$P^\sadj$ is computed
formally from $p_i^\sadj=i\partial_{p_i}$ and
$\partial_{p_i}^\sadj=p_i/i$; in particular
$(p_i\partial_{p_i})^\sadj=p_i\partial_{p_i}$ \cite{Macdonald95}.
This makes the adjunction map an involution as well as an algebra
anti-automorphism of~$W_p$.  Note that, although adjunction extends
to~$W_p(t)$ by setting~$t_i^\sadj=t_i$, no adjoint for
the~$\partial_{t_i}$ can be defined in any consistent way.  Assume
that an adjoint~$\partial_{t_i}^\sadj$ existed.  For reasons to be
explained later, this adjoint has to be of the form
$\alpha\partial_{t_i}+\beta t_i+\gamma$ for complex constants
$\alpha$, $\beta$, $\gamma$, with $\alpha\beta\neq0$.  Now, for any
series $F$ and~$G$ we have $\rsc{\partial_{t_i}\cdot F}G=\rsc
F{\partial_{t_i}^\sadj\cdot G}$.  Choose any non-zero series~$F$
independent of~$t_i$; then by the method of variation of parameters for
series, one finds a series~$G$ satisfying $\partial_{t_i}^\sadj\cdot
G=F$.  Upon evaluation, we obtain $0=\rsc FF\neq0$, a contradiction.

We now proceed to outline the algorithm for the simple case, meaning that
from this point on we elect to have $F\in\field[[p]]$, i.e.,
$F$~independent of~$t$.  The condition on~$F$ that it does not
involve~$t$ implies that $\partial_{t_i}\cdot F=0$ for $i$ from~1
to~$k$.  We can use this fact to simplify our calculations.  In this
case, we consider a different annihilator, $\ann_{W_p}F$, hereafter
denoted~$J_F$.  Note that $J_F=\cI_F\cap W_p$.

This allows us to determine the action of combinations of $P\in
J_F^\sadj$ and $Q\in\cI_G$.  For example, given any $S\in W_p$, $T\in
W_{p,t}(t)$, and $U\in W_t(t)$,
\[
\rsc F{(P^\sadj SU+TQ)\cdot G}=
\rsc{S^\sadj P\cdot F}{U\cdot G}+\rsc F{TQ\cdot G}
=0.
\]
It follows that, if we can find a combination such that
$\sum_jP_j^\sadj S_jU_j+\sum_jT_jQ_j=R\in W_t$, we have $0=\rsc
F{R\cdot G}=R\cdot\rsc FG$.  Note that each~$P_j^\sadj S_j$ is an
element of~$J_F^\sadj$ while each~$T_jQ_j$ is an element of~$\cI_G$.
Therefore, we conduct our search for an element of $\ann_{W_t}\rsc FG$
by determining a non-zero element of $\bigl(J_F^\sadj
W_t(t)+\cI_G\bigr)\cap W_t$.  We shall prove in
Section~\ref{sec:proof1and3} that such an element exists.  Basically,
the goal of our algorithms is to compute sufficiently many non-zero
elements of $\bigl(J_F^\sadj W_t(t)+\cI_G\bigr)\cap W_t$ so as to
generate a D-finite description of the scalar product.

Note, however, that while~$\cI_G$ is a left $W_{p,t}(t)$ ideal,
$J_F^\sadj W_t(t)$ is a {\em right\/} $W_{p,t}(t)$-ideal and the
sums~$P+Q$ for $P\in J_F^\sadj W_t(t)$ and $Q\in\cI_G$ do not form an
ideal.  This problem is very similar to the problem of creative
telescoping: given an ideal~$\cI\subset W_{p,t}(t)$, the aim in the
first step of this method is to determine an element of
$\partial_pW_{p,t}(t)+\cI$ that does not involve~$p$. There also,
$\partial_pW_{p,t}(t):=\sum_j\partial_{p_j}W_{p,t}(t)$ is a right
ideal. The algorithm we present thus bears a non-fortuitous
resemblance with that of~\cite{Takayama90b}: in this reference,
truncations of the left ideal~$\cI$ and of the right
ideal~$\partial_pW_{p,t}(t)$ at a given total degree in
$p,\partial_p,\partial_t$ are recombined linearly, this for higher and
higher truncation degrees until the corresponding truncation of the
intersection $\bigl(\partial_pW_{p,t}(t)+\cI\bigr)\cap W_t$ is
non-trivial.  In our situation, we determine truncations of the left
ideal~$\cI_G$ and the right ideal~$J_F^\sadj W_t(t)$ at a given
truncation order, recombine those two vector spaces linearly, and
iterate over higher and higher truncation orders until the
corresponding truncation of $\bigl(J_F^\sadj W_t(t)+\cI_G\bigr)\cap
W_t$ is a D-finite description.

To some extent, the approach of the present paper also shares features
with that in~\cite{OaTa99}.  However, this reference focuses on
determining a bound on a truncation order that permits to compute
generators of an intersection $L=\bigl(\partial_pW_{p,t}+I\bigr)\cap
W_t$ for a given ideal~$I$ of~$W_{p,t}$, and also generators for a
whole free resolution of~$L$.  From there, the cohomology groups of
the module-theoretic integral~$W_t/L$ of the quotient
module~$W_{p,t}/I$ are derived.  Roughly speaking, we are not
concerned here with more than the first cohomology group, and
furthermore, we treat the similar but different problem for ideals
of~$W_{p,t}(t)$ and intersections in~$W_t(t)$.

Being a module over~$W_t(t)$, the sum~$J_F^\sadj W_t(t)+\cI_G$ is a
vector space over~$\field(t)$.  It is this second structure that is
adapted to our method.  We could try using the module structure in
this section, but this would not generalize to the case when also~$F$
depends on~$t$.  The idea is to use $\field(t)$-linear algebra in the
vector space structure to eliminate the $\partial_{p_i}$ and
$p_i$. Roughly speaking, we incrementally generate lines in a matrix
corresponding to generators of~$J_F^\sadj W_t(t)+\cI_G$, and perform
Gaussian elimination to remove the monomials involving $p$
and~$\partial_p$.

The main loop of the algorithm considers monomials of increasing
degree with respect to any ordering on the monomials
in~$p,\partial_p,\partial_t$.  We use the notation~$\preceq$ to denote
the monomial comparison associated with this ordering.  We reduce each
monomial~$\alpha$ with respect to (the Gr\"obner bases for)
$\cI_F^\sadj$ and~$\cI_G$.  Note that the chosen monomial ordering is
the same for both $\cI_G$ and~$\cI_F^\sadj$.  Equivalently,
the remainder of the reduction of a monomial~$\alpha$
with respect to~$\cI_F^\sadj$ can be viewed as the adjoint of the
remainder of the reduction of~$\alpha^\sadj$ with respect to~$\cI_F$.
However, to reflect the fact that adjunction modifies the variables,
when reducing with respect to~$\cI_F$ we need to use a different
order, specifically, the ordering~$\preceq_\sadj$ defined by
$\beta_1\preceq_\sadj\beta_2$ on~$W_p$ if and only if
$\beta_1^\sadj\preceq\beta_2^\sadj$.  In our implementation, we use
the ordering $\operatorname{DegRevLex}(\partial_p>p>\partial_t)$ which
sorts by total degree first, breaking ties by a reverse lexicographic
order on the variables.  The order~$\preceq_\sadj$ is then
$\operatorname{DegRevLex}(p>\partial_p)$.

Once we have computed these values, we add two rows (and for
sufficiently large~$\alpha$ only one column) in a matrix where we
perform Gaussian elimination to cancel entries corresponding to
monomials involving $p$ and~$\partial_p$.

We now state the algorithm more formally as Algorithm~\ref{thm:algo1},
followed by an example in the next section. After this example, we
describe the modifications necessary to handle specific cases more
efficiently, and how to treat the general case. The proofs that these
algorithms work and terminate are delayed until Section~\ref{sec:proofs}.

\begin{algo}[Scalar Product]\label{thm:algo1}
\mbox{}\\
{\sc Input:} Symmetric functions\/ $F\in\field[[p]]$ and\/
$G\in\field[t][[p]]$, both D-finite in\/~$p,t$, given by D-finite
descriptions in\/~$W_p$ and\/~$W_{p,t}(t)$, respectively.\\
{\sc Output:} A D-finite description of\/ $\rsc FG$ in\/~$W_t(t)$.
\begin{enumerate}

\item Determine a Gr\"obner basis\/~$\cG_G$ for the left ideal\/
$\ann_{W_{p,t}(t)}G$ with respect to any monomial ordering\/~$\preceq$,
as well as a Gr\"obner basis\/~$\cG_{F^\sadj}$ for the right ideal\/
$\ann_{W_p}F^\sadj$ with respect to the monomial ordering induced
by\/~$\preceq$ on\/~$W_p$;
    \label{item:gb-comput-in-algo-1}

\item $B:=\{\}$;

\item Iterate through each monomial\/~$\alpha$ in\/
$p,\partial_p,\partial_t$;

  \begin{enumerate}
  \item Write\/ $\alpha=\beta\gamma$ with\/ $\beta\in W_p$ and\/
  $\gamma\in\field[\partial_t]$;
  \item $\alpha_F:=\bigl(\beta-(\beta\redu_\preceq\cG_{F^\sadj})\bigr)\gamma$;
  \item $\alpha_G:=\alpha-(\alpha\redu_\preceq\cG_G)$;
  \item Introduce\/ $\alpha_F$ and\/ $\alpha_G$ as two new elements into\/ $B$
  and reduce so as to eliminate\/~$p,\partial_p$;
	\label{item:after-intro-in-algo-1}
  \item  Compute the dimension of the ideal generated by\/ $B\cap
  W_t(t)$. If this dimension is~0, break and output\/ $B\cap W_t(t)$.
	\label{item:breaking-condition-in-algo-1}
  \end{enumerate}

\end{enumerate}
\end{algo}

Notice, if~$m=1$, as is the case in our examples, there is only one
variable~$t$, and the dimension condition
in~\eqref{item:breaking-condition-in-algo-1} is simplified to:
\begin{center}\em
If\/ $B$~contains a non-zero element\/~$P$ from\/~$W_t(t)$, break
and return\/~$P$.
\end{center}

Note that Step~(\ref{item:gb-comput-in-algo-1}) requires to determine
both ideals $\ann_{W_{p,t}(t)}G$ and~$\ann_{W_p}F$, not just
$\ann_{W_{p,t}(p,t)}G$ and $\ann_{W_p(p)}F$.  In other words, one
generally needs to pass from a D-finite description~$\{P_i\}$
generating the ideal $\ann_{W_p(p)}F$ as $\sum_iW_p(p)P_i$ to a
set~$\{Q_i\}$ generating the ideal $\ann_{W_p}F=W_p\cap\ann_{W_p(p)}F$
as $\sum_iW_pQ_i$, and similarly for~$G$.  The operation of computing
such intersections is called {\em Weyl closure}, in the terminology
of~\cite{Tsai00,Tsai02}.  It is a non-obvious task, owing to the
change of module structure (coefficients in~$W_p(p)$ are replaced with
coefficients in~$W_p$).  Algorithms are provided
in~\cite{Tsai00,Tsai02}.

Sometimes, the input set~$\{P_i\}$ already constitutes a generating
set for the Weyl closure.  In this case, one can skip
Step~(\ref{item:gb-comput-in-algo-1}) of the algorithm.  This is the
case in our examples.

The remainder of the reduction with respect to the Gr\"obner
basis~$\cG_G$ is a multivariate analogue of the remainder of the
Euclidean division.  It is such that for any~$\alpha$,
$\alpha_{\cG}=\alpha-(\alpha\redu\cG)$ belongs to the ideal generated
by~$\cG$.  A~similar statement holds for~$\cG_F$.

For this description we have assumed that Gr\"obner bases could be
computed for both left and right ideals.  If they can only be computed
on one side, say for left ideals, then the operators~$\alpha_F$ can be
obtained as follows: first, determine the monomial
ordering~$\preceq_\sadj$ induced by adjunction on~$W_p$ viewed as a
left structure from the ordering~$\preceq$ on~$W_p$ viewed as a right
structure; then, replace the Gr\"obner basis~$\cG_{F^\sadj}$ with the
Gr\"obner basis~$\cG_F$ for the left ideal~$\ann_{W_p}F$ with respect
to~$\preceq_\sadj$; $\alpha_F$~is then computed as
$\bigl(\beta-(\beta^\sadj\redu_{\preceq_\sadj}\cG_F^\sadj)\bigr)\gamma$.
This way we get $\cG_{F^\sadj}=(\cG_F)^\sadj$.

We represent the basis~$B$ as a matrix, with columns indexed by
monomials in the~$p_i$'s, the~$\partial_{p_i}$'s, and
the~$\partial_{t_i}$'s.  Each row in the matrix corresponds to the row
vector of the coefficients of some element of~$B$ with regard to the
indexing monomial basis.  Introducing an element into the basis
consists of adding a new row at the bottom of the matrix, performing
row reduction (also known as Gaussian elimination), and then returning
the new matrix as the updated basis.  In practice, $B$~can be handled
(not inefficiently) by a Gr\"obner basis computation with respect to a
monomial ordering that eliminates the~$p_i$'s and
the~$\partial_{p_i}$'s, performing calculations in the free
$\field[t]$-module with a basis the list of indexing monomials.

Finally, some remembering can be done at Step~(3b) to avoid reducing
the same~$\beta$ again and again, for different~$\alpha$'s involving
the same~$\beta$.

\section{Example: $k$-Regular Graphs}\label{sec:example}

The enumeration of regular graphs has been treated by a number of
authors \cite{Comtet74,Gessel90,GoJaRe83,ReWo80}.  We present it here
because of its expository value and as it is the simplest in a family
of examples.  After expressing the problem as a scalar product, we
describe in detail how our algorithm treats it.  We conclude this
section with an indication of how the scenario may be generalized.

\subsection{A generating series for graphs as a scalar product}

Recall from the introduction that a generating series for the set of
all finite simple graphs labeled with integers from
$\bN\setminus\{0\}$ is
\[
G(x)=\sum_{G\in\cG}\prod_{(i,j)\in E(G)}x_ix_j=\prod_{i<j}(1+x_ix_j),
\]
under the encoding that a graph on $n$~vertices $i_1,\dots,i_n$ of
respective valencies $v_1,\dots,v_n$ contributes a
monomial~$x_{i_1}^{v_1}\dots x_{i_n}^{v_n}$.  We can similarly make a
generating function for graphs with multiple edges (multigraphs) by
\[
M(x)=\prod_{i<j}\frac1{(1-x_ix_j)},
\]
for an edge~$(i,j)$ of a graph with multiplicity~$m$ contributes a
monomial~$x_i^mx_j^m$ and any non-negative multiplicity is allowed.

Clearly both $G$ and~$M$ are symmetric functions, and in fact, we have
the relations $G=e[e_2]$ and~$M=h[e_2]$, as determined by a method
that we discuss in Section~\ref{sec:species}.  Both are easily
rewritten in terms of the~$p_i$'s:
\begin{equation}\label{eqn:G}
G=\exp\left(\sum_i(-1)^i(p_i^2-p_{2i})/2i\right)
\quad\text{and}\quad
M=\exp\left(\sum_i\left(p_i^2+p_{2i}\right)/2i\right).
\end{equation}

In any given term, the degree of~$x_k$ gives the valency of
vertex~$k$.  So, for example, the coefficient~$g_n$ of~$x_1\dotsm x_n$
in~$G$, hereafter denoted~$[x_1\dotsm x_n]G$, gives the number of
1-regular graphs, or perfect matchings on the complete graph on
$n$~vertices, and in general the coefficient $g^{[k]}_n=[x_1^k\dotsm
x_n^k]G$, also given as $[m_{k^n}]G$, gives the number of $k$-regular
graphs on $n$~vertices.  By virtue of Eq.~\eqref{eq:scalhm},
coefficient extraction amounts to a scalar product, and the generating
function~$G_k(t)$ of $k$-regular graphs is given by
\begin{multline}\label{eq:kreg}
G_k(t):=\sum_ng^{[k]}_n\frac{t^n}{n!}=\rsc G{H_k},
\qquad\text{where}\\
H_k(t):=\sum_n h_{k^n}\frac{t^n}{n!}=\sum_n\frac{(h_kt)^n}{n!}=\exp(h_kt).
\end{multline}

Now, since $h_k=\sum_{\lambda\vdash k} p_\lambda/z_\lambda$ (where the
sum is over all partitions~$\lambda$ of~$k$), the exponential
generating function~$H_k(t)$ is also $\exp\bigl(t\sum_{\lambda\vdash
n}p_\lambda/z_\lambda\bigr)$, an exponential in a finite number of
$p_i$'s. By Property~(3) in Theorem~\ref{thm:dprops}, this is
D-finite.  Further, as a result of scalar product
property~\eqref{eq:scalp}, we can rewrite Eq.~\eqref{eq:kreg} as
\begin{equation}\label{eq:kregp}
G_k(t)=\rsc
{\exp\left(
   \sum_{i\text{ even},\ i\le k}(-1)^{i/2}\frac{p_i^2}{2i}+\frac{p_i}i
  +\sum_{i\text{ odd},\ i\le k}{\frac{p_i^2}{2i}}
\right)\!}
{\,\exp\left(t\sum_{\lambda\vdash k}\frac{p_\lambda}{z_\lambda}\right)}
\end{equation}
and now by Theorem~\ref{thm:rsc_pres_df} this generating function
$G_k(t)$ is D-finite.

Note how the closed form for~$G$ in~(\ref{eqn:G}), in infinitely many
variables, and the closed form for~$H_k(t)$ in~(\ref{eq:kreg}), in
terms of the~$h$'s, have led to the scalar product~(\ref{eq:kregp})
between two closed forms, explicitly written in terms of finitely
many~$p_i$ for each~$k$.  This reduction is what has made the
algorithm applicable.

\subsection{Effective Computation for $k=2$}
To illustrate a typical calculation, we calculate $G_2(t)$, the
generating function for 2-regular graphs which, according to
Eq.~\eqref{eq:kregp}, is determined by
\[
G_2(t)=\rsc{\exp\bigl((p_1^2-p_2)/2-p_2^2/4\bigr)}{\exp\bigl(t(p_1^2+p_2)/2\bigr)}.
\]
Algorithm~1 calculates that $G_2(t)$ satisfies the differential
equation\[2(1-t)G'_2(t)-t^2G_2(t)=0,\] which is easily solved to find
$G_2(t)=e^{-\frac14t(t+2)}/\sqrt{1-t}$.

In order to appeal to Algorithm~1, set $F=\exp((p_1^2-p_2)/2-p_2^2/4)$
and $G=\exp(t(p_1^2+p_2)/2)$ and determine the Gr\"obner bases $\cG_F$
and $\cG_G$ of their annihilating ideals respectively:
\[\cG_F=
  \{
     \underline{p_2}+2\partial_{p_2}+1,
     \underline{p_1}-\partial_{p_1}
  \}
\quad\text{and}\quad
\cG_G=\{
        2\underline{\partial_{p_2}}-t,
        \underline{\partial_{p_1}}-tp_1,
        \underline{p_1^2}+p_2-2\partial_t\},
\]
where $\cG_F$~is a Gr\"obner basis with respect to the degree reverse
lexicographical monomial ordering such that
$p_1>p_2>\partial_{p_1}>\partial_{p_2}$ and $\cG_G$~is a Gr\"obner
basis with respect to the degree reverse lexicographical monomial
ordering such that $\partial_{p_1}>\partial_{p_2}>p_1>p_2>\partial_t$.
(Leading monomials with respect to the monomial ordering are
underlined.)  Before proceeding, the set~$\cG_F$ is converted by
adjunction into a Gr\"obner basis~$\cG_F^\sadj$ with respect to the
degree reverse lexicographical monomial ordering such that
$\partial_{p_1}>\partial_{p_2}>p_1>p_2$:
\[\cG_F^\sadj=
  \{
     2\underline{\partial_{p_2}}+p_2+1,
     \underline{\partial_{p_1}}-p_1
  \}.
\]
(The reader should not get confused by the peculiar situation of this
example: here, adjunction has not changed the polynomials, except for
signs, but this is only a coincidence.)

The initial value of $B$ is the empty set.  For the sake of the
example, we shall iterate on monomials~$\alpha$ according to the
degree reverse lexicographical order such that
$\partial_t>\partial_{p_2}>p_2>\partial_{p_1}>p_1$, and perform
reductions when inserting into the basis according to the elimination
order sorting first by the degree reverse lexicographical order such
that $\partial_{p_2}>p_2>\partial_{p_1}>p_1$, and breaking ties by the
degree in~$\partial_t$.

We now briefly sketch the run of the algorithm until $\alpha$~becomes
$p_1\partial_{p_1}$ and then illustrate the steps of the main loop in
more details.

For $\alpha=1$ and~$\alpha=p_1$, the algorithm inserts no polynomial
into the basis~$B$.  The next iteration of the loop,
for~$\alpha=\partial_{p_1}$, produces
$\alpha_F=\underline{\partial_{p_1}}-p_1$, which is inserted into~$B$
as is, and $\alpha_G=\underline{\partial_{p_1}}-tp_1$, whose insertion
puts~$p_1$ into~$B$.  Next, the case $\alpha=p_2$ inserts no
polynomial before, for~$\alpha=\partial_{p_2}$,
$\alpha_F=2\underline{\partial_{p_2}}+p_2+1$ gets inserted as is, and
the insertion of~$\alpha_G=2\underline{\partial_{p_2}}-t$
puts~$\underline{p_2}+(t+1)$ into~$B$.  The iteration
for~$\alpha=\partial_t$ has no effect on~$B$.  For~$\alpha=p_1^2$,
$\alpha_F=0$~is not inserted, and
$\alpha_G=\underline{p_1^2}+p_2-2\partial_t$ gets inserted in the form
$\underline{p_1^2}-2\partial_t-(t+1)$.

At this point, the algorithm is about to
treat~$\alpha=p_1\partial_{p_1}$ and the value of~$B$ is
\begin{equation}\label{eq:B}
B=
\left\{
        \underline{\partial_{p_1}}-p_1,
        \underline{p_1},
        2\underline{\partial_{p_2}}+p_2+1,
        \underline{p_2}+(t+1),
        \underline{p_1^2}-2\partial_t-(t+1)
\right\},
\end{equation}
where we have written elements in the order of introduction into the
set.  In matrix notation, the column vector of elements of~$B$ reads:
\begin{equation*}
\begin{pmatrix}
0&0&0&1&-1&0&0\\
0&0&0&0&1&0&0\\
0&2&1&0&0&0&1\\
0&0&1&0&0&0&t+1\\
1&0&0&0&0&-2&-(t+1)
\end{pmatrix}
\begin{pmatrix}
p_1^2\\
\partial_{p_2}\\
p_2\\
\partial_{p_1}\\
p_1\\
\partial_t\\
1
\end{pmatrix}
\end{equation*}
Here, we have chosen to keep the rows in the order of creation by the
algorithm and to sort the column according to the monomial order used
by the elimination step.  Observe that in this way, no two rows have
their left-most non-zero entry on the same column: simply reordering
rows would put the matrix in row echelon form.

Then, the algorithm computes
\begin{multline*}
\alpha_F=\alpha-(\alpha\redu_\preceq\cG_F^\sadj)=\alpha-(\alpha^\sadj\redu_{\preceq_\sadj}\cG_F)^\sadj
=\underline{p_1\partial_{p_1}}-p_1^2+1\\
\text{and}\qquad
\alpha_G=\alpha-(\alpha_{\redu_\preceq\cG_G})= \underline{p_1\partial_{p_1}}+tp_2-2t\partial_t.
\end{multline*}
(Note that $\alpha_F$~is really~$(\partial_{p_1}-p_1)p_1$, an element
of the \emph{right\/} ideal generated by~$\cG_F^\sadj$.)  Next, we
update~$B$ to include these two values. We insert~$\alpha_F$ into~$B$
after one reduction, leading to
\[B:=B\cup\{\underline{p_1\partial_{p_1}}-2\partial_t-t\}.\]
In matrix notation, this insertion adds a new column to the left of
the matrix, corresponding to the new monomial~$p_1\partial_{p_1}$, and
one more row at the bottom of the matrix,
$(\begin{smallmatrix}1&0&0&0&0&0&-2&-t\end{smallmatrix})$.  Then the
algorithm inserts~$\alpha_G$. Its leading monomial~$p_1\partial_{p_1}$
is already present in~$B$, leading to an initial reduction
to~$t\underline{p_2}+2(1-t)\partial_t+t$.  One final reduction by
$t$~times the pre-last element in Eq.~\eqref{eq:B} results in the step
\[B:=B\cup\{2(1-t)\underline{\partial_t}-t^2\}.\]
The intersection of this and~$W_t(t)$ is non-trivial, and the algorithm
outputs $2(1-t)\partial_t-t^2$. We conclude that $G_2(t)$ satisfies
the differential equation
\[2(1-t)G_2'(t)-t^2G_2(t)=0.\]
Table~\ref{table-kreg} summarizes the results by the same algorithm for $k=2,3,4$.
These match with the results in~\cite{GoJaRe83}.

\begin{table}
\begin{center}
\begin{tabular}{rc}
\hline
&2-regular graphs\\\cline{2-2}
$\phi_0$&$-t^2$\\
$\phi_1$&$-2t+2$\\
$\phi_2$&0\\
\hline
&3-regular graphs\\\cline{2-2}
$\phi_0$&$t^3(t^4+2t^2-2)^2$\\
$\phi_1$&$-3(t^{10}+6t^8+3t^6-6t^4-26t^2+8)$\\
$\phi_2$&$-9t^3(t^4+2t^2-2)$\\
\hline
&4-regular graphs\\\cline{2-2}
$\phi_0$&$-t^4(t^5+2t^4+2t^2+8t-4)^2$\\
$\phi_1$&\parbox{12cm}{$-4(t^{13}+4t^{12}-16t^{10}-10t^9-36t^8-220t^7-348t^6$\\\hbox{}\hfill$-48t^5+200t^4-336t^3-240t^2+416t-96)$}\\
$\phi_2$&$16t^2(t-1)^2(t^5+2t^4+2t^2+8t-4)(t+2)^2$\\
\hline
\end{tabular}
\bigskip
\caption{Differential equation $\phi_2 G_k''+\phi_1 G_k'+\phi_0 G_k=0$
satisfied by $G_k(t)$, $k 
= 2,3,4$.\label{table-kreg}}
\end{center}
\end{table}

\subsection{Efficient enumeration of $k$-regular graphs}

An efficient procedure for the enumeration of $k$-regular graphs is
immediately derived from the differential equations for the generating
series of $k$-regular graphs collected in Table~\ref{table-kreg}.
Indeed, one simply needs to convert the differential equation
for~$G_k(t)$ into a recurrence relation for its
coefficients~$g^{[k]}_n$ and to determine sufficiently many starting
values $g^{[k]}_0$, $g^{[k]}_1$,~\dots\@{} Then, one can efficiently
compute~$g^{[k]}_n$ for any~$n$ by unrolling the recurrence.

Implementations are available to help with this approach.  For
example, the Maple package {\tt gfun}\footnote{This package is part of
the {\tt algolib} library, which is available at {\tt
http://\discretionary{}{}{}algo\discretionary{}{}{}.inria\discretionary{}{}{}.fr/\discretionary{}{}{}packages/}.}
by Salvy and Zimmerman~\cite{SaZi94} contains commands dedicated to
the conversion step and the iterative calculations based on a linear
recurrence.  Computations in the case~$k=4$ result in a recurrence
relation of order~15 already published by Read and Wormald
\cite{ReWo80} and can be found as a formula accompanying sequence
number A005815 in Sloane's encyclopedia of integer sequences
\cite{Sloane03}.  From this recurrence relation and initial terms, it
is then a matter of seconds to compute the exact integer values for
hundreds of terms in the sequence.

It should be stressed that this method proves much more efficient than
the direct computation of the scalar product based on a termwise
expansion and application of formula~\eqref{eq:scalp}.  For example,
Stembridge's implementation in the package SF for symmetric function
manipulation in Maple~\cite{Stembridge95} already requires several
minutes to compute the~$g^{[4]}_n$ for~$n$ up to~15, and becomes
unsuitable to handle the symmetric functions that would be necessary
to obtain~$g^{[4]}_{20}$.  Far from showing any weakness of SF's
general approach, this illustrates the computational progress provided
by our techniques in the specific setting of differentiably finite
series.

\subsection{Generalization}\label{sec:species}

The series given by Eq.~\eqref{eqn:G} is determined combinatorially in
a direct fashion using the theory of species~\cite{BeLaLe98}. This can
be extended naturally to handle a wider family of combinatorial
structures, such as hypergraphs, set covers, latin rectangles. For an
in-depth treatment, consult~\cite{Mishna03}.

\section{Hammond Series}\label{sec:hammond} 

In the example above, it turned out that except for monomials of
degree~1, we needed only examine the two monomials $p_1^2$
and~$p_1\partial_{p_1}$ in order to reach the solution. However,
depending on the monomial ordering, the algorithm might well consider
many monomials before it adds the ones that eliminate the $p_i$'s
and~$\partial_{p_i}$'s. The problem becomes far more serious as the
number of variables and the degree of the monomials increase. It turns
out that in the common case when the scalar product is of the
type~$\rsc F{H_k(t)}$ it is possible to modify the approach and
eliminate the~$p_i$ and the~$\partial_{p_i}$ in a more efficient
manner using the {\em Hammond series\/}\footnote{In
\cite[Sec.~3.5]{GoJa83} this is referred to as the {\em Gamma
series\/} of~$F$.} (or H-series) introduced by Goulden, Jackson, and
Reilly in~\cite{GoJaRe83}: for $F\in\field[[p_1,p_2,\dotsc]]$, the Hammond
series of~$F$ is defined as
\[
\cH(F)(t_1, t_2,\dotsc)=\rsc F{\sum_\lambda h_\lambda t^\lambda/m(\lambda)!},
\]
where the sum is over all partitions, and if $\lambda=1^{m_1}\dotsm
k^{m_k}$ then $t^\lambda=t_1^{m_1}\dotsm t_k^{m_k}$ and
$m(\lambda)!=m_1!\,m_2!\dotsm m_k!$.  These are very closely related to
the Hammond operators, defined by Hammond~\cite{Hammond83} and used
extensively by MacMahon~\cite{MacMahon60}.  A Hammond operator can be
described as~$h_\lambda^\sadj$, and thus the Hammond series of~$F$
with all of the $t$~variables set to~1 results essentially in a sum of
Hammond operators acting on~$F$.

Observe that the generating function for $k$-regular graphs is
\[G_k(t)=\cH(G)(0,\dots, 0, t, 0, \dotsc)\]
where the~$t$ occurs in position~$k$. This is true for any generating
function which takes the form $\rsc F{H_k(t)}$ for some~$F$.

A theorem from~\cite{GoJaRe83} is specially useful: Goulden, Jackson,
and Reilly's H-series theorem states that $\cH(\partial_{p_n}\cdot F)$
and~$\cH({p_n}F)$ can be expressed in terms of
the~$\partial_{t_i}\cdot\cH(F)$'s. In terms of Gr\"obner bases, this
corresponds to introducing the additional variables $t_1$,~\dots,
$t_k$ (instead of~$t=t_k$ alone) and work with the series
$\cH_k(t_1,\dots,t_k)=\sum_\lambda h_\lambda t^\lambda/m(\lambda)!$
with sum over partitions~$\lambda$ whose largest part is~$k$ (instead
of working with the univariate~$H_k(t)$). The H-series theorem
therefore implies that for an appropriate monomial order, there is a
Gr\"obner basis of the ideal~$I_{\cH_k}$ of all operators of $W_{p,t}$
annihilating~$\cH_k$, with elements of the form
\begin{equation}\label{gbhk}
p_i-P_i(t,\partial_t),\quad
\partial_{p_i}-Q_i(t,\partial_t),\qquad i=1,\dots,k,
\end{equation}
where all the $P_i$ and~$Q_i$ are polynomials in~$t,\partial_t$.

The algorithm in this case is as follows.
\begin{algo}[Hammond Series]\label{thm:algo_jgr} 
\mbox{}\\
{\sc Input:} An integer~$k$, and $F\in\field[[p_1,\dots,p_n]]$.\\
{\sc Output:} A differential equation satisfied by
\[\rsc F{\sum_i h_{k^i}t_k^i}=\cH(F)(0,\dots,0,t_k,0,\dotsc)\]
where $t_k$~is in position~$k$.
\begin{enumerate}
\item Compute $\cG_F$, a Gr\"obner basis for the left ideal~$J_F$
annihilating~$F$ in~$W_p$; 
\item Compute $\cG_{\cH_k}$, a Gr\"obner basis of the
form~\eqref{gbhk};
\item For each~$U\in\cG_F$, compute $r_U\in W_t$ as the reduction
of~$U^\sadj$ by~$\cG_{\cH_k}$ for an order which
eliminates~$p,\partial_p$.  Let~$R_0$ be the set of~$r_U$'s;
\item For~$i$ from 1 to $k-1$ 
eliminate~$\partial_{t_i}$ from~$R_{i-1}$ and set~$t_i=0$ in
the resulting polynomials; call~$R_i$ the new set;
\item Return~$R_{k-1}$.
\end{enumerate}
\end{algo}

As with Algorithm~1, the first step is to determine an annihilating
ideal in~$W_p$.  Again, one can possibly first determine a D-finite
description and use Weyl closure~\cite{Tsai00,Tsai02} to obtain the
annihilating ideal.

After Step~(3), all the~$p_i$'s and~$\partial_{p_i}$'s have been
eliminated and $R_0$~contains a set of generators of a D-finite
$W_t(t)$-ideal annihilating~$\rsc F{\cH_k}$. Then, in order to obtain
differential equations satisfied by the specialization
at~$t_1=\dots=t_{k-1}=0$, Step~(4) proceeds in order by eliminating
differentiation with respect to~$t_i$ and then setting~$t_i=0$ in the
remaining operators.

Note that the Gr\"obner basis of Step~(2) can be precomputed for the
required~$k$'s (although most of the time is actually spent in Step~(4)).

In order to compute the elimination in Step~(4), one should not
compute a Gr\"obner basis for an elimination order, since this would
in particular perform the unnecessary computation of a Gr\"obner basis
of the eliminated ideal. Instead, one can modify the main loop in the
Gr\"obner basis computation so that it stops as soon as sufficient
elimination has been performed or revert to skew elimination by the
non-commutative version of the extended Euclidean algorithm as
described in~\cite{ChSa98}.  This is the method we have adopted in the
example session given in Appendix~B\footnote{An implementation of the
algorithms presented here is available in the Maple package
ScalarProduct available at {\tt
http://\discretionary{}{}{}algo.\discretionary{}{}{}inria.\discretionary{}{}{}uqam.\discretionary{}{}{}fr/\discretionary{}{}{}{
}mishna}.}.

This calculation is comparatively rapid since the size of the basis is
greatly reduced. Further, the basis grows smaller as the algorithm
progresses, on account of setting variables to~0. We can compute the
case of 4-regular graphs in a second, instead of a couple of minutes
using the general algorithm.  The 5-regular expression requires
significantly more computation time, and we could not compute it.

A mathematically equivalent but slightly faster way of performing
Step~(3) is to compute~$r_U$ by simply replacing each monomial
$p_1^{\alpha_1}\dotsm
p_n^{\alpha_n}\partial_{p_1}^{\beta_1}\dotsm\partial_{p_n}^{\beta_n}$
in~$U$ with the product $Q_n^{\beta_n}\dotsm
Q_1^{\beta_1}P_n^{\alpha_n}\dotsm P_1^{\alpha_1}$.

In order to explain the relative speed of Algorithm~2, compared to
Algorithm~1, it needs to be said that the Hammond series approach searches
a smaller space, which can well result in a differential equation of
order higher than that obtained by Algorithm~1.  This occurs, for
instance, in the case of 4-regular graphs: Algorithm~2 returns a
differential equation of order~3 only when that returned by
Algorithm~1 is of order~2.

In the same vein, note that the order in which the eliminations are
done in Step~(4) could be changed, possibly leading to a different
(but correct) output.

\subsection{Proof of Termination and Correctness}

Termination of Algorithm~2 is obvious.  On the other hand, the full
proof of correctness requires technical results to be proved in
Section~\ref{sec:proofs}.  The following corollary articulates a
property of D-finite functions in the simple language of symmetric
functions and D-finite descriptions, and is a corollary of
Proposition~\ref{thm:struct-S'} that will be proved independently.

\begin{cor}\label{thm:IFplusIG}
Let $F$ and~$G$ be D-finite symmetric series in
$\field[[p_1,\dots,p_n]]$ and
$\field[t_1,\dots,t_k][[p_1,\dots,p_n]]$, respectively, with
corresponding annihilators $J_F\subset W_p$ and $\cI_G\subset
W_{p,t}(p,t)$.  Under these conditions, the vector space
\[\left(J_F^\sadj W_t(t)+\cI_G\right)\cap W_t(t)\]
is non-trivial and contains a D-finite description of $\rsc FG$.
\end{cor}

\begin{prop}
Algorithm~2 terminates and is correct. 
\end{prop}

\begin{proof}
First, we remark that for fixed~$k$,
\[\cH_k(t_1,\dots,t_k)=\exp\left(\sum_{j=1}^kh_jt_j\right)\]
is a D-finite symmetric series by Theorem~\ref{thm:dprops} since
each~$h_j$ is a finite combination of $p_1$,~\dots, $p_n$. Thus,
$f=H(F)(t_1,\dots,t_k)=\rsc{\cH_k(t_1,\dots,t_k)}F$ is a D-finite
function of $t_1$,~\dots, $t_k$, by Theorem~\ref{thm:rsc_pres_df}.

We proceed by proving the following invariant of the main loop: the
set~$R_{i-1}$ generates a D-finite description of
$\cH(F)(0,\dots,0,t_i,t_{i+1},\dots,t_k)$.  This establishes the
result since it implies that $R_{k-1}$ contains a D-finite description
of $\cH(F)(0,\dots,0,t_k)$, in this case, a single differential
equation. This is precisely what the algorithm claims to determine.

To prove the base case of this invariant, note that $R_0$ contains the
generators of the intersection $\left(J_F^\sadj
W_t(t)+\cI_{\cH_k}\right)\cap W_t(t)$.  We appeal to
Corollary~\ref{thm:IFplusIG}, to conclude that $R_0$~contains a
D-finite description of $\cH(F)(t_1,\dots,t_k)$.

The general case is proven with the known result~\cite{ChSa98} that
given a D-finite description of a function $F(x_1, x_2,\dots,x_n)$,
one can compute the D-finite description of $F(x_1,\dots,x_{n-1},0)$,
for example, by first eliminating~$\partial_{x_n}$, removing factors
of~$x_n$ in the remaining polynomials, and finally, setting~$x_n=0$ in
the equations, precisely the process outlined in Algorithm~2.
\end{proof}

\section{Example: $k$-Uniform Tableaux}
\label{sec:young}

Another family of combinatorial objects whose generating function can
be resolved with our method is a certain class of Young tableaux,
namely $k$-uniform Young tableaux.

For a partition~$\lambda=(\lambda_1,\dots,\lambda_k)\vdash n$, a Young
tableau of shape~$\lambda$ is an array $T=(T_{i,j})$ of positive
integers~$T_{i,j}$ defined when $1\leq i\leq k$ and~$1\leq
j\leq\lambda_i$.  When a Young tableau is strictly increasing on each
of its rows and columns ($T_{i,j}<T_{i+1,j}$ and~$T_{i,j}<T_{i,j+1}$,
whenever this makes sense) and the $n$~integers~$T_{i,j}$ are all
integers from~1 to~$n$, it is called standard.

Standard Young tableaux are in direct correspondence with many
different combinatorial objects. For example, Stanley~\cite{Stanley99}
has studied the link between standard tableaux and paths in Young's
lattice, the lattice of partitions ordered by inclusion of
diagrams. This link was generalized by Gessel~\cite{Gessel93b} to
tableaux with repeated entries. Gessel remarks that such paths have
arisen in the work of Sundaram on the combinatorics of representations
of symplectic groups~\cite{Sundaram90}.

The weight of a tableau is $\mu=(\mu_1,\dots, \mu_k)$ where $\mu_1$ is
the number of 1's, $\mu_2$ is the number of 2's,~etc., in the tableau
entries.  Here we consider Young tableaux that are column strictly
increasing and row weakly increasing, and with weight
$\mu=1^k2^k\dotsm n^k$: each entry appears $k$~times.  We call Young
tableaux with these properties {\em$k$-uniform}.  These correspond to
paths in Young's lattice with steps of length~$k$. The set of
$k$-uniform tableaux of size~$kn$ is also in bijection with symmetric
$n\times n$ matrices with non-negative integer entries with each row
sum equal to~$k$.  Gessel notes that for fixed~$k$, the generating
series of the number of $k$-uniform tableaux is
D-finite~\cite{Gessel90}. Our method makes this effective.

Two observations from~\cite{Macdonald95} are essential. First,
$[x_1^{\mu_1}\dotsm x_k^{\mu_k}]s_\lambda$ is the number of (column
strictly increasing, row weakly increasing) tableaux with weight
$\mu$. Secondly,
\[
\sum_\lambda s_\lambda = h[e_1+e_2]
=\exp\left(\sum_i p_i^2/2i + \sum_{i\text{ odd}}p_i/i\right),
\]
which is D-finite.  Define $y_n^{[k]}$ to be the number of $k$-uniform
tableaux of size $kn$, and let~$Y_k$ be the generating series of these
numbers.  The previous two observations imply
\begin{equation}\label{eq:kunif}
Y_k(t)=\sum_n y_n^{[k]} t^k
=\rsc{{\exp\left(\sum_{i=1}^k p_i^2/2i +
\sum_{i\text{ odd}}^kp_i/i\right)}}{\sum_n h_{k^n} t^n},
\end{equation}
This problem is well-suited to our methods since again we treat an
exponential of a polynomial in the $p_i$'s, with an explicit closed
form in terms of~$k$ for this polynomial.

Calculating the equations for $k=1,2,3,4$ is fast with either
Algorithm~1 or Algorithm~2. The resulting differential equations are
listed in Table~\ref{ftab}.  For $k=1,2$ these results agree with
known results~\cite{Gupta68,Stanley99}, and are the entries A000085
and A000985 respectively in Sloane's encyclopedia of integer sequences
\cite{Sloane03}.  The first few values of $y_n^{[k]}$ are summarized
in the following table. For $k=3,4$ these appear to be new.

Concerning the dual problem, where instead $n$~is fixed and
$k$~varies, the sequences $\bigl(y_n^{[k]}\bigr)_k$ appear
respectively as A019298, A053493, and A053494 for $n=3,4,5$.  Stanley
\cite[Prop.~4.6.21]{Stanley86} reports that the generating functions
$G_n(x)=\sum_k y_n^{[k]} x^k$ are rational with denominator of the
form $(1-x)^s(1-x^2)^t$ where $s$ and~$t$ are positive integers.

\begin{table}
\center
\begin{tabular}{rc}
\hline
&1-uniform tableaux\\\cline{2-2}
$\phi_0$&$-(t-1)$\\
$\phi_1$&$1$\\
$\phi_2$&0\\
\hline
&2-uniform tableaux\\\cline{2-2}
$\phi_0$&$t^2(t-2)$\\
$\phi_1$&$-2(t-1)^2$\\
$\phi_2$&0\\
\hline
&3-uniform tableaux\\\cline{2-2}
$\phi_0$&$(t^{11}+t^{10}-6t^9-4t^8+11t^7-15t^6+8t^5-2t^3+12t^2-24t-24)$\\
$\phi_1$&$-3t(t^{10}-2t^8+2t^6-6t^5+8t^4+2t^3+8t^2+16t-8)$\\
$\phi_2$&$9t^3(-t^2-2+t+t^4)$\\
\hline
&4-uniform tableaux\\\cline{2-2}
$\phi_i$&(See Appendix~\ref{sec:fourreg})\\
\hline
\end{tabular}
\bigskip
\caption{Differential equation $\phi_2Y_k''+\phi_1Y_k'+\phi_0Y_k=0$
satisfied by $Y_k(t)$, $k=1,\dots,4$.\label{ftab}} 
\end{table}

\begin{table}
\begin{small}
\begin{center}
\begin{tabular}{ll}
\hline
$k$& $y_0^{[k]}, y_1^{[k]}, y_2^{[k]},\dotsc$\\
\hline
1&1, 1, 2, 4, 10, 26, 76, 232, 764, 2620, 9496, 35696, 140152, 568504\smallskip\\
2&\parbox{12cm}{1, 1, 3, 11, 56, 348, 2578, 22054, 213798, 2313638, 27627434, 360646314,\\\hbox{}\hskip1em5107177312, 77954299144\medskip}\\
3&\parbox{12cm}{1, 1, 4, 23, 214, 2698, 44288, 902962, 22262244, 648446612, 21940389584,\\\hbox{}\hskip1em849992734124\medskip}\\
4&\parbox{12cm}{1, 1, 5, 42, 641, 14751, 478711, 20758650, 1158207312, 80758709676,\\\hbox{}\hskip1em6877184737416, 701994697409136\smallskip}\\
\hline
\end{tabular}
\bigskip
\caption{The number, $y_n^{[k]}$, of $k$-uniform tableaux of size $kn$.}
\end{center}
\end{small}
\end{table}

\section{Algorithm for Scalar Product: the General Situation}\label{sec:general}

So far, we have limited the scope of the algorithms to pairs of
D-finite symmetric functions where only one of the two functions
depends on the variables $t_1$,~\dots, $t_k$.  While this is
sufficient in many applications, it is possible to modify Algorithm~1
in order to accommodate the~$t_i$'s in both functions and thus make
the full power of Theorem~\ref{thm:rsc_pres_df} effective.  While no
additional ideas are to be used, the description of the algorithm is
more technical.

Algorithm~1 manipulates monomials $\alpha$ and reduces them modulo the
ideals $\cI_F$ and~$\cI_G$ in order to determine equations of the
form
\begin{equation}\label{eqn:alg1}
\rsc F{\bigl(\alpha-(\alpha\redu_\preceq\cI_F^\sadj)\bigr)\cdot G}=0
\quad\text{and}\quad
\rsc F{\bigl(\alpha-(\alpha\redu_\preceq\cI_G)\bigr)\cdot G}=0,
\end{equation}
where on the left, $\alpha$~supposedly does not involve any of
the~$\partial_{t_i}$'s.  What makes the situation of Algorithm~1 and
the left-hand identity in~\eqref{eqn:alg1} simple is the assumption
that $F$~does not depend on~$t$, making the action of~$W_t$ on~$\rsc
FG$ act on the right-hand argument only.  The difficulty in
generalizing lies in that now, the action of~$\partial_{t_i}$ on~$F$
may be non-trivial and must be considered in the differentiation rule
for scalar products,
\begin{equation}\label{eq:dif}
\partial_{t_i}\cdot\rsc FG=\rsc{\partial_{t_i}\cdot F}G+\rsc
F{\partial_{t_i}\cdot G},
\end{equation}
which itself stems from the differentiation rule for usual products on
the level of coefficients.

The idea is therefore to manipulate operators in {\em three\/} sets
of~$\partial_{t_i}$'s: one which acts on the full scalar product $\rsc
FG$, and one for each of its components, acting directly on the
component. To facilitate the description of this situation, we denote
the former by~$\partial_{t_i}$, the one acting on the left component
by~$\partial_{\ell_i}$, and the one acting on the right
component~$\partial_{r_i}$.  Using this notation, we wish to view
Eq.~\eqref{eq:dif} as
\begin{equation}\label{eqn:dif3}
\partial_{t_i}=\partial_{\ell_i}+\partial_{r_i}.
\end{equation}

We thus modify Algorithm~1 by enlarging the family of monomials over
which we iterate, and use Eq.~\eqref{eqn:dif3} to eliminate
the~$\partial_{\ell_i}$'s before we begin Gaussian elimination.  Here, we
iterate over monomials $\alpha\partial_\ell^\beta\partial_r^\gamma$ of
the free commutative monoid $[p,\partial_p,\partial_\ell,\partial_r]$
with~$\alpha\in[p,\partial_p]$ to examine the following
generalizations of~Eq.~\eqref{eqn:alg1}:
\begin{multline}\label{eqn:generalized-alg1}
\rsc{\bigl(\alpha^\sadj\partial_t^\beta-(\alpha^\sadj\partial_t^\beta\redu\cG_F)\bigr)\cdot
F}{\partial_t^\gamma\cdot G}=0\\
\text{and}\qquad
\rsc{\partial_t^\beta\cdot
F}{\bigl(\alpha\partial_t^\gamma-(\alpha\partial_t^\gamma\redu\cG_G)\bigl)\cdot
G}=0,
\end{multline}
or, with a change of notation,
\begin{multline*}
\bigl(\alpha^\sadj\partial_\ell^\beta-(\alpha^\sadj\partial_\ell^\beta\redu\cG_F)\bigr)\partial_r^\gamma\cdot\rsc FG=0\\
\text{and}\qquad
\partial_\ell^\beta\bigl(\alpha\partial_r^\gamma-(\alpha\partial_r^\gamma\redu\cG_G)\bigl)\cdot\rsc FG=0.
\end{multline*}
Upon making use of Eq.~\eqref{eqn:dif3} and applying adjunction to the
first equation in Eq.~\eqref{eqn:generalized-alg1}, we get a linear
combination of terms of the form $\partial_t^{\beta'}\cdot\rsc
F{\alpha'\cdot G} $ with coefficients in~$\field[t]$, where
$\beta'\in\bN^k$, and $\alpha'\in W_{p,t}(t)$. The algorithm proceeds
as before by performing Gaussian elimination over~$\field(t)$ to
eliminate $p,
\partial_p$, and $\partial_r$.  In our implementation, the monomial
order~$\preceq$ is
$\operatorname{DegRevLex}(\partial_r>\partial_\ell>\partial_p>p)$.  The
method is summarized in Algorithm~\ref{thm:algo2}.

\begin{algo}[General Scalar Product]\label{thm:algo2}
\mbox{}\\
{\sc Input:} $F\in\field[t][[p]] $ and\/ $G\in\field[t][[p]]$, both
D-finite in\/~$p,t$, given by D-finite descriptions
in\/~$W_{p,t}(t)$.\\
{\sc Output:} A D-finite description of\/ $\rsc FG$ in\/~$W_t(t)$.
\begin{enumerate}

\item Determine a Gr\"obner basis\/~$\cG_G$ for the left ideal\/
$\ann_{W_{p,t}(t)}G$ with respect to any monomial ordering\/~$\preceq$,
as well as a Gr\"obner basis~$\cG_{F^\sadj}$ for the right ideal\/
$\ann_{W_{p,t}}F^\sadj$ with respect to the same ordering;

\item $B:=\{\}$;

\item Iterate through each monomial\/~$\alpha$ in\/
$p,\partial_p,\partial_\ell,\partial_r$ in any order;

  \begin{enumerate}
        \item $\alpha_l:=\alpha|_{\partial_\ell=\partial_t,\partial_r=1}$;
        \item $\alpha_F:=\alpha_l-(\alpha_l\redu_\preceq\cG_{F^\sadj})$;
        \item $\alpha_r:=\alpha|_{\partial_r=\partial_t,\partial_\ell=1}$;
        \item $\alpha_G:=\alpha_r-(\alpha_r\redu_\preceq\cG_G)$;
        \item Introduce\/
	$(\alpha_F|_{\partial_\ell=\partial_t-\partial_r})
	(\alpha|_{p=\partial_p=\partial_\ell=1})$ and\/
	$(\alpha|_{p=\partial_p=\partial_r=1})\alpha_G$ into\/~$B$
	and reduce so as to eliminate\/ $p,\partial_p,\partial_r$;
		\label{item:after-intro-in-algo-2}
        \item Compute the dimension of the ideal generated by\/
        $B\cap W_t(t)$. If this dimension is~$0$, break and output\/
        $B\cap W_t(t)$.
		\label{item:breaking-condition-in-algo-2}
  \end{enumerate}

\end{enumerate}
\end{algo}

As in Algorithm~1, if~$m=1$, there is only one variable~$t$, and the
condition in~\eqref{item:breaking-condition-in-algo-2} is simplified
to:
\begin{center}\em
If\/ $B$~contains a non-zero element\/~$P$ from\/~$W_t(t)$, break and
return\/~$P$.
\end{center}

The same remarks as those made after Algorithm~1 at the end of
Section~\ref{sec:alg} also apply here.

\section{Termination and Correctness of Algorithms 1 and~3}\label{sec:proofs}

\subsection{Sketch of the proof}\label{sec:proof1and3}
The common goal of Algorithms 1 and~3 is to find differential
equations satisfied by $\rsc FG$, which is equivalent to non-zero
elements in~$W_t$ which annihilate $\rsc FG$.  Although Algorithm~1 is
a specialization of Algorithm~3, parts of the proof would become
artificially more involved if restricted to the simple case.  We thus
treat both algorithms simultaneously.  The discussion at the beginning
of Section~\ref{sec:alg} has illustrated how to manipulate the
annihilators of $F$ and~$G$ to determine a combination $P^\sadj
S+TQ\in W_t$ with $P\in\cI_F^\sadj$, $Q\in\cI_G$, $S\in W_p(t)$, $T\in
W_{p,t}(t)$, which annihilates~$\rsc FG$.  Not all of the elements in
$\ann_{W_t}\rsc FG$ are of this form, however, as the following simple
example illustrates. If $F=p_1-p_2$ and $G=p_1+p_2/2$, then $\rsc
FG=1-1=0$ and thus $1\in\ann_{W_t}\rsc FG$. However, it can be
established that 1 can not be written as a combination of the form
$P^\sadj S+T Q$ for those $F$ and~$G$. Nonetheless, we show that the
annihilating elements that can be written this way form a non-trivial
subideal of $\ann_{W_t}\rsc FG$, which we generate with the
algorithms.

Although the problem of finding differential equations appears at first
inherently analytic in nature, we rephrase it algebraically into a
question amenable to the theory of D-modules. The adjunction
properties of the scalar product are naturally accommodated by tensor
products.  Specifically, the proof below centers around a certain
$W_t$-module~$S$ whose elements are tensors, and where, for example,
\[
(i^{-1}p_i\cdot u)\otimes v=(u\cdot\partial_{p_i})\otimes
v=u\otimes(\partial_{p_i}\cdot v),
\]
which corresponds to the equivalence $\rsc{(i^{-1}p_i)\cdot F}G =\rsc
F{\partial_i\cdot G}$.  (See also
Eq.~(\ref{eqn:modified-Laplace}--\ref{eqn:tensor-product-rules-end})
below.)  On the other hand, the $\partial_{\ell_i}$ and~$\partial_{r_i}$
that are involved in the description of Algorithm~3 really are the
operators $\partial_{t_i}\otimes1$ and~$1\otimes\partial_{t_i}$ acting
on~$S$, respectively, where 1's~denote identity maps.
 
The module $S$ can be expressed in terms of the ideal
$\ann_{W_t}(F^\sadj\otimes G)$, itself contained in $\ann_{W_t}\rsc
FG$.  The former ideal is non-trivial and in fact, is sufficient to
describe the scalar product as holonomic, a property whose definition
is recalled shortly and which implies D-finiteness.  In fact, we
show that the algorithms calculate a Gr\"obner basis for
$\annFoGt$, in other words a D-finite description of the scalar
product~$\rsc FG$.

The main result  is summarized by the following theorem.
\begin{thm}\label{thm:main}
Suppose $F$ and~$G$ are symmetric functions subject to the conditions
of Algorithm~1 (resp.\ Algorithm~3).  Then, Algorithm~1 (resp.\
Algorithm~3) determines, in finite time, a Gr\"obner basis for a
non-zero D-finite ideal contained in\/ $\ann_{W_t(t)}\rsc FG$.
\end{thm}

The notion of holonomy to be used in the proof follows
\cite{Borel87,Coutinho95}.  Introduce a filtration of~$W_t$ by the
$\field$-vector spaces~$F_d$ of all operators in~$W_t$ of total degree
at most~$d$ in~$t,\partial_t$.  These spaces are finite-dimensional,
of dimension $\binom{d+2k}{2k}=O\bigl(d^{2k}\bigr)$ as $d$~tends to
infinity.  A $W_t$-module $M=\sum_iW_t\cdot g_i$ generated by a finite
family of generators~$g_i$ is holonomic whenever the $\field$-vector
spaces $\sum_iF_d\cdot g_i$ have dimension growing like
$O\bigl(d^k\bigr)$.  A function of~$t$ that is an element of a
holonomic $W_t$-module is called holonomic.  From the definition, it
is a basic result that a holonomic function is D-finite; the converse
is a more difficult result to be found in \cite[Th.~2.4 and
Appendix~6]{Takayama92}.  Similar definitions apply to
$W_{p,t}$-modules, with a dimension growth of $O\bigl(d^{k+n}\bigr)$
in place of $O\bigl(d^k\bigr)$.

The discussion so far has not relied on the definition of the scalar
product. Rather, remark that Algorithms 1 and~3 are essentially
parameterized by the adjunction property of the scalar product of
symmetric functions, and can easily be redefined and adapted to other
adjunctions.  It suits our needs for the proof to consider adjoints for
the usual scalar product of functions, $\rsp fg:=\int f(x)g(x)\,dx$. To
avoid confusion, we notationally distinguish $\rsp fg$ from~$\rsc FG$ for
the two scalar products, as well as $\uadj$ from~$\sadj$ for the
respective adjunction operations.

Indeed, guided by existing results concerning the preservation of
holonomy under operations involving the usual scalar product, we link
the symmetric case to the usual one with a map from one adjunction to
the other.  This reduction also demonstrates how algorithms analogous
to Algorithms 1 and~3 for other scalar products could be shown to
terminate with the correct output.  (See
Section~\ref{sec:other-scalar-products}.)
 
To make this comparison more intuitive, we could identify $\rsc FG$
with the integral
\[
\int_{\bR^n}\cL\bigl(q\mapsto F(q_1,2q_2,\dots,nq_n)\bigr)(p)
\,G(p)\,dp_1\dotsm dp_n,
\]
where $\cL$ is the modified Laplace transform
\[\cL(F)(p)=\int_{\bR^n}F(q)e^{-(p_1q_1+\dots+p_nq_n)}\,dq,\]
which satisfies
\[
\cL\bigl(q\mapsto
q_iF(q)\bigr)(p)=-(\partial_{p_i}\circ\cL)(F)(p).
\]
Notice,
for example:
\begin{multline}\label{eqn:modified-Laplace}
\rsc{i^{-1}p_i\cdot F}G=
\int_{\bR^n}\cL\bigl(q\mapsto q_iF(q_1,\dots,nq_n)\bigr)(p)
\,G(p)\,dp_1\dotsm dp_n\\
=-\int_{\bR^n}(\partial_{p_i}\circ\cL)(F)(p)
\,(\partial_{q_i}\cdot G)(p)\,dp_1\dotsm dp_n\\
=\int_{\bR^n}\cL\bigl(q\mapsto F(q_1,\dots,nq_n)\bigr)(p)
\,(\partial_{q_i}\cdot G)(p)\,dp_1\dotsm dp_n
=\rsc F{\partial_{p_i}\cdot G}.
\end{multline}
Formally, we must work on the level of abstract modules, however.
This avoids situations where the integral is not convergent or the
Laplace transform is not defined as a function.

Thus, to prove Theorem~\ref{thm:main}, we show
Corollary~\ref{thm:elements_of} below which states that
$\ann_{W_t}\left(F^\sadj\otimes G\right)$ is a non-zero subideal of
$\ann_{W_t}\rsc FG$ such that the quotient
$W_t/\ann_{W_t}\left(F^\sadj\otimes G\right)$ is a holonomic module.
This is done in several stages. First, in Section~\ref{sec:reduction},
we define~$S$, the algebraic structure in which our calculations take
place, and prove that it is holonomic by reducing the problem to the
usual scalar product analogue, where similar results are known.  This
analogue is detailed in Section~\ref{sec:usual}.  Next, in
Section~\ref{sec:quotient} we express $S$ as a
quotient. Corollary~\ref{thm:elements_of} follows from this
discussion.  Finally, to conclude that the algorithm terminates, we
relate~$S$ to the algorithm in more detail and prove in
Section~\ref{sec:termination} that all of the generators are
determined in finite time. Together, these results prove
Theorem~\ref{thm:main} and thus the correctness and termination of
Algorithms 1 and~3.

\subsection{The scalar product of symmetric 
functions}\label{sec:reduction}

We now formally define the $W_t$-module~$S$. Begin with
$U=W_{p,t}\cdot F$ and $V=W_{p,t}\cdot G$, two holonomic
$W_{p,t}$-modules.  We shall denote by $U^\sadj$ the adjoint module
of~$U$: as $\field$-vector spaces, $U=U^\sadj$, and a right
$W_p[t]$-action is defined on~$U^\sadj$ by $u\cdot P=P^\sadj\cdot u$
for any $u\in U^\sadj$ and $P\in W_p[t]$, where the last operation is
taken for the left structure of~$U$.  Set $S$ as the tensor product
$U^\sadj\otimes_{W_p[t]}V$, which makes it a $\field[t]$-module.  This
has the desirable effect of encoding the scalar product adjunction
relations: for all $u\in U$ and all $v\in V$,
\begin{gather}
\label{eqn:tensor-product-rules-begin}
(\partial_{p_i}\cdot u)\otimes v
=(u\cdot\partial_{p_i}^\sadj)\otimes v
=(u\cdot i^{-1}p_i)\otimes v =u\otimes(i^{-1}p_i\cdot v),\\
(p_i\cdot u)\otimes v
=(u\cdot p_i^\sadj)\otimes v
=(u\cdot i\partial_{p_i})\otimes v
=u\otimes(i\partial_{p_i}\cdot v),\\
\label{eqn:tensor-product-rules-end}
t_i\cdot(u\otimes v)
=(t_i\cdot u)\otimes v
=(u\cdot t_i)\otimes v
=u\otimes(t_i\cdot v).
\end{gather}

To endow $S$ with a $W_t$-module structure, let $\partial_{t_i}$ act
on a pure tensor $u\otimes v$ by
\begin{equation}\label{eq:dt-action}
\partial_{t_i}\cdot(u\otimes v)=(\partial_{t_i}\cdot u)\otimes
v+u\otimes(\partial_{t_i}\cdot v),
\end{equation}
and extend to~$S$ by $\field$-linearity.  In other words,
$\partial_{t_i}=\partial_{\ell_i}+\partial_{r_i}$ after defining
$\partial_{\ell_i}=\partial_{t_i}\otimes1$ and
$\partial_{r_i}=1\otimes\partial_{t_i}$, where 1's~are identity maps.

Armed with this definition and Theorem~\ref{thm:T_hol} (formally
stated and proven independently in Section~\ref{sec:usual}), we prove
that $S$~is holonomic. Theorem~\ref{thm:T_hol} is an analogous result
for the usual scalar product, corresponding adjunction, and
corresponding adjoint module~$M^\uadj$ of a module~$M$.  It states
that for holonomic $M$ and~$N$, $M^\uadj\otimes_{W_p[t]}N$ is a
holonomic $W_t$-module under the action of~$\partial_{t_i}$ given
by~\eqref{eq:dt-action}.  We shall appeal to this theorem with an
appropriate choice for $M$ and~$N$.

To determine the relationship between the two scalar products and make
our choice for $M$ and~$N$, we compare both adjunction operations.  In
the symmetric case, adjunction is defined as the
anti-automorphism~$\sadj$ which maps $p_i$ to $i\partial_{p_i}$ and
$\partial_{p_i}$ to $i^{-1}p_i$, for all $i$, and the usual scalar
product adjunction is defined as the anti-automorphism~$\uadj $ which
maps $\partial_{p_i}$ to $-\partial_{p_i}$, and leaves the $p_i$
variables unchanged.  One way to connect both adjunctions is to factor
$\sadj$ into the composition of three algebra morphisms:

\begin{enumerate}
\item the automorphism~$\tau$ mapping $(p_i,\partial_i)$
to~$(ip_i,i^{-1}\partial_i)$.  This corresponds to the dilation which
maps a function~$F$ to~$p\mapsto F(p_1,2p_2,\dots,np_n)$;

\item the automorphism~$\cF$ mapping $(p_i,\partial_i)$
to~$(-\partial_i,p_i)$ and named `Fourier transform' in D-module
theory (see \cite[proof of Th.~3.1.8]{Borel87} or
\cite[p.~39]{Coutinho95}).  Informally speaking, this corresponds to
mapping a function~$F$ to its Laplace transform~$\cL(F)$;

\item the anti-automorphism~$\uadj $ mapping $(p_i,\partial_i)$
to~$(p_i,-\partial_i)$. 
\end{enumerate}
The important property to note is that each of these three maps
preserves holonomy since they preserve total degree, hence are
filtration-preserving bijections.  A direct calculation on $p_i$
and~$\partial_i$ verifies that $\sadj=\uadj\circ\cF\circ\tau$, so that
the composite~$\sadj$ also is a holonomy-preserving linear bijection.
Thus, we introduce two holonomic modules, $M=(\cF\circ\tau)(U)$ also
denoted $U^{\cF\circ\tau}$, and~$N=V$, so as to appeal to
Theorem~\ref{thm:T_hol}.  One concludes that
\begin{equation}\label{eq:adj-fact}
S=U^\sadj\otimes_{W_p[t]}V
=\left(U^{\cF\circ\tau}\right)^\uadj\otimes_{W_p[t]}V
=M^\uadj\otimes_{W_p[t]}N
\end{equation}
is a holonomic $W_t$-module. After we have described the quotient
structure of $S$ in Section~\ref{sec:quotient}, this information will be
used to prove that $\annFoG$~is non-trivial and that the quotient
module $W_t/\annFoG$ is holonomic, a fact we use to show that the
algorithms terminate.

\subsection{Preservation of holonomy under the usual scalar product}\label{sec:usual}
In the previous section, we reduced the proof of the holonomy of
$S=U^\sadj\otimes_{W_p[t]}V$ to an analogous result in terms of the
usual scalar product, to be proven in this section: the module
$T=M^\uadj\otimes_{W_p[t]}N$ is holonomic when $M$ and~$N$ are.

The following notion will be used in the proof: the integral of a
$W_{p,t}$-module~$P$, denoted $\int P=\int P\,dp_1\dotsm dp_n$, is
defined as $P\bigm/\bigl(\sum_i\partial_{p_i}\cdot P\bigr)$.  It is
the image of composed maps: the Fourier transform~$\cF$, the inverse
image~$\pi_*$ under the projection~$\pi$ from~$\field^{n+m}$
to~$\field^n$ defined by $\pi(p,t)=t$, and the inverse Fourier
transform.  Specifically we have, $\int P=\cF^{-1}\pi_*\cF(P)$.  These
maps preserve holonomy (see
\cite[Th.~3.3.4]{Borel87} or
\cite[Th.~18.2.2 and Sec.~20.3]{Coutinho95}), so that the integral of
a holonomic $W_{p,t}$-module is a holonomic $W_t$-module.  (See also
\cite[Th.~3.1.8]{Borel87}.)

The module~$T$ fits naturally in between an existing
holonomy-preserving surjection from the $W_t$-module $\int
M\otimes_{\field[p,t]}N$ to the space $\rsp MN$.  Factoring this map
to pass through $T=M^\uadj\otimes_{W_p[t]}N$ yields:

\begin{equation}\label{eqn:reg_link}
\int M\otimes_{\field[p,t]}N
\stackrel{\phi}{\longrightarrow}
M^\uadj\otimes_{W_p[t]}N
\stackrel{\psi}{\longrightarrow}
\rsp MN,
\end{equation}
where $\psi$~surjectively maps~$m\otimes n$ to~$\rsp mn$, and
$\phi$~is a natural $W_t$-linear surjection that we are about to
define in the course of the next theorem.  After proving that the
first module in~\eqref{eqn:reg_link} is holonomic, the surjectivity
of~$\phi$ implies the holonomy of~$T$.

\begin{thm}\label{thm:T_hol} Suppose that $M$ and~$N$ are two holonomic
$W_{p,t}$-modules, and define~$T$ as $M^\uadj\otimes_{W_p[t]}N$.
Then, $T$ is a holonomic $W_t$-module under the action
of~$\partial_{t_i}$ given by
\[
\partial_{t_i}\cdot(m\otimes n)=(\partial_{t_i}\cdot m)\otimes n +
m\otimes(\partial_{t_i}\cdot n).
\]
\end{thm}

\begin{proof}
First, we focus our attention on the module $\int
M\otimes_{\field[p,t]}N$ in~\eqref{eqn:reg_link}.  Consider the
$W_{p,t}$-module $P:=M\otimes_{\field[p,t]}N$, with action
of~$\partial_{p_i}$ defined by $\partial_{p_i}\cdot(m\otimes
n)=(\partial_{p_i}\cdot m)\otimes n+m\otimes(\partial_{p_i}\cdot n)$,
and action of~$\partial_{t_i}$ defined similarly. We can also write
this as the inverse image $\iota^*\left(M\otimes_{\field}N\right)$,
where $\iota$ is the map from $\field^{m+n}$ to $\field^{(n+m)+(n+m)}$
which sends~$(p,t)$ to~$(p,t,p,t)$.  The advantage of the second
presentation is that the holonomy of~$P$ is obtained from the
holonomic closure under inverse image under embeddings (see
\cite[Th.~3.2.3]{Borel87} or
\cite[Sec.~15.3 and Ex.~15.4.5]{Coutinho95}) and the holonomic closure
under tensor product over~$\field$
\cite[Cor.~13.4.2]{Coutinho95}.  Therefore, $\int P$~is also holonomic.

Next, we define a $W_t$-linear surjection to $T$.  Define a map
from~$M\times N$ to~$T$ which sends $(m,n)$ to $m\otimes n$.  This map
is $\field[p,t]$-balanced, $\field[p,t]$-bilinear, and surjective. By
the universality of the tensor product, this induces a surjective
map~$\phi$ from~$P=M\otimes_{\field[p,t]}N$ to~$T$.  Observe that each
derivation~$\partial_{p_i}$ maps~$P$ into the kernel of~$\phi$, as the
following calculation indicates:
\begin{multline*}
\phi\bigl(\partial_{p_i}\cdot(m\otimes n)\bigr)
    =\phi\bigl((\partial_{p_i}\cdot m)\otimes n+m\otimes(\partial_{p_i}\cdot n)\bigr)\\
    =(\partial_{p_i}\cdot m)\otimes n+m\otimes(\partial_{p_i}\cdot n)
    =m\otimes(-\partial_{p_i}\cdot n)+m\otimes(\partial_{p_i}\cdot n)
    =0.
\end{multline*}
In other words, $\sum_i\partial_{p_i}\cdot P\subset\ker\phi$, and thus
$\phi$~also induces a well-defined surjective map from~$\int P$
to~$T$.  Any good filtration of~$\int P$ will induce a good filtration
for~$T$ (see \cite[Prop.~1.11]{Borel87} or
\cite[Lemma~7.5.1]{Coutinho95}). Thus, $T$ is finitely generated with
dimension bounded by that of $\int P$.  Therefore, $T$~is holonomic.
\end{proof}

\subsection{The quotient structure of $S$}\label{sec:quotient}
Subsequent developments to express~$S$ as a quotient involve modules
over~$W_{p,t}$ and ideals of~$W_{p,t}$, rather than~$W_{p,t}(t)$.  We
therefore introduce the annihilators $I_F=\ann_{W_{p,t}}F$ and
$I_G=\ann_{W_{p,t}}G$, to be used in place of
$\cI_F=\ann_{W_{p,t}(t)}F$ and $\cI_G=\ann_{W_{p,t}(t)}G$,
respectively.  Note that $I_F=\cI_F\cap W_{p,t}$ and
$\cI_F=\field(t)\otimes_{\field[t]}I_F$, and similarly for~$G$.
Finally, although adjunction has not been defined for~$\partial_t$, we
use the notation~$W_{p,t}^\sadj$ to denote $W_{p,t}$ endowed with both
a structure of $W_t$-module on the left and a structure of
$W_p[t]$-module on~the~right.

\begin{prop}
The module~$S=(W_{p,t}\cdot F)^\sadj\otimes_{W_p[t]}(W_{p,t}\cdot G)$
is isomorphic to
\[(\WoW)/(\FplusG).\]
\end{prop}
\begin{proof}
The $W_t$-module $S=U^\sadj\otimes_{W_p[t]}V$ is also a
$W_{p,t}^\sadj\otimes_{W_{p[t]}}W_{p,t}$-module.  As such, it is
generated by~$F^\sadj\otimes G$.  Consider the two exact sequences of
respectively right and left $W_p[t]$-modules
\[
\begin{array}{ccccccccc}
\label{ex_seq_F}
0&\rightarrow&I_F^\sadj &\xrightarrow{\rho}&W_{p,t}^\sadj &
  \xrightarrow{\alpha}&U^\sadj&\rightarrow&0,
\\
0&\rightarrow&I_G&\xrightarrow{\sigma}&W_{p,t}&
  \xrightarrow{\beta}&V&\rightarrow&0,
\end{array}
\]
where $\alpha(P)=F^\sadj\cdot P$, $\beta(Q)=Q\cdot G$, and $\rho$
and~$\sigma$ are inclusions.  (Here, $F$ and~$F^\sadj$ denote the same
element of the set~$U$, but we write~$F^\sadj$ when viewed as an
element of the right module~$U^\sadj$, $F$~when viewed as in the left
module~$U$.)  We combine them to make a third exact sequence:
\begin{equation}\label{eqn:ses}
\begin{array}{cccccccc}
\ker(\alpha\otimes\beta) &\rightarrow   &\WoW
&\xrightarrow{\alpha\otimes\beta}       &S&\rightarrow& 0,\\
&&P\otimes Q             &\longmapsto   &(F^\sadj\cdot P)\otimes(Q\cdot G)
\end{array}
\end{equation}
where, by \cite[II.59, Proposition~6]{Bourbaki70},
\[
\ker(\alpha\otimes\beta)
=\im(\rho\otimes 1_{W_{p,t}})+\im(1_{W_{p,t}^\sadj}\otimes\sigma)
=\FplusG
\] 
as $\field[t]$-modules. We conclude that, as $W_t$-modules,
\begin{multline*}
S\simeq(\WoW)/\ker(\alpha\otimes\beta)\\
\simeq(\WoW)/(\FplusG).
\end{multline*}
\end{proof}

To be more explicit, note that this isomorphism maps the class
of~$1\otimes1$ in the quotient to~$F^\sadj\otimes G\in S$.  Remark
also that, as $W_t$-modules,
\begin{multline*}
\ker(\alpha\otimes\beta)
=\bigl\{P\otimes Q\in W_{p,t}^\sadj\otimes W_{p,t}:(\alpha\otimes\beta)(P\otimes Q)=0\bigr\}\\
=\bigl\{P\otimes Q\in W_{p,t}^\sadj\otimes W_{p,t}:(F^\sadj\cdot P)\otimes(Q\cdot G)=0\bigr\}\\
=\bigl\{P\otimes Q\in W_{p,t}^\sadj\otimes W_{p,t}:(P\otimes Q)\cdot(F^\sadj\otimes G)=0\bigr\}\\
=\annFxG,
\end{multline*}
so that we also have
\begin{equation}\label{eq:3-formulations}
\annFxG=\ker(\alpha\otimes\beta)=\FplusG.
\end{equation}

\begin{prop}\label{thm:struct-S'}
The $W_t$-module~$S'=W_t\cdot(F^\sadj\otimes G)$ is a submodule
of~$S$, isomorphic to
\[\Wt\bigm/\bigl((\FplusG)\cap\Wt\bigr),\]
where $\Wt\simeq W_t$~is the smallest $\field$-subalgebra
of~$W_{p,t}^\sadj\otimes_{W_p[t]}W_{p,t}$ generated by $\field[t]$,
$1\otimes\partial_{t_1}+\partial_{t_1}\otimes1$,~\dots,
$1\otimes\partial_{t_k}+\partial_{t_k}\otimes1$.  In the simplified
situation when $I_F=\partial_tW_{p,t}+W_tJ_F$ for~$J_F=\ann_{W_p}F$,
$S'$~is isomorphic~to
\[W_t\bigm/\bigl((W_tJ_F^\sadj+I_G)\cap W_t\bigr).\]
\end{prop}

We first prove this proposition, then in the next section we discuss
how to connect the description of~$S'$ above directly to the algorithm
and how to apply it to show that the algorithms terminate.

\begin{proof}
The annihilator of~$F^\sadj\otimes G$ in~$W'_t\cdot(F^\sadj\otimes G)$
\[\ann_{\Wt}(F^\sadj\otimes G)=\ann_{\WoW}(F^\sadj\otimes G)\cap\Wt.\]
In view of the action of~$W_t$ on~$S'$ through the isomorphism between
$W_t$ and~$\Wt$, we thus have that $S'$~is isomorphic to
$W_t/\ann_{W_t}(F^\sadj\otimes G)$, itself isomorphic~to
\[
\Wt/\ann_{\Wt}(F^\sadj\otimes G)
=\Wt/\bigl(\ann_{\WoW}(F^\sadj\otimes G)\cap\Wt\bigr).
\]
Owing to~\eqref{eq:3-formulations}, this proves the general quotient
expression for~$S'$ in the proposition statement.

Now, to prove the formula in the simpler case, observe that
when~$I_F=\partial_t W_{p,t}+W_tJ_F$,
\begin{multline*}
I_F^\sadj\otimes_{W_{p[t]}}W_{p,t}
=\partial_tW_{p,t}^\sadj\otimes_{W_{p[t]}}W_{p,t}+
W_tJ_F^\sadj\otimes_{W_{p[t]}}W_{p,t}\\
=\partial_tW_t\otimes_{\field[t]}W_{p,t}+
W_t\otimes_{\field[t]}W_tJ_F^\sadj
\end{multline*}
while $W_{p,t}^\sadj\otimes_{W_p[t]}I_G=W_t\otimes_{\field[t]}I_G$,
whence the relation
$\ker(\alpha\otimes\beta)=\partial_tW_t\otimes_{\field[t]}W_{p,t}+
W_t\otimes_{\field[t]}(W_tJ_F^\sadj+I_G)$.
Since~$\WoW=W_t\otimes_{\field[t]}W_{p,t}$, we obtain
\[S\simeq W_{p,t}/(W_tJ_F^\sadj+I_G),\]
as $(W_t\otimes_{\field[t]}W_{p,t})/\ker(\alpha\otimes\beta)
\simeq(\field[t]\otimes_{\field[t]}W_{p,t})
/\bigl(\field[t]\otimes_{\field[t]}(W_tJ_F^\sadj+I_G)\bigr)
\simeq W_{p,t}/(W_tJ_F^\sadj+I_G)$.
Following these isomorphisms, $\Wt$~can be identified as the copy
of~$W_t$ included in~$W_{p,t}$ in the last quotient above.  Therefore,
the submodule~$S'$ of~$S$ is isomorphic to the quotient announced in
the proposition statement.
\end{proof}

\begin{cor}\label{thm:elements_of}
The ideal\/ $\annFoG$ is: 
\begin{enumerate}
\item isomorphic to $(\FplusG)\cap\Wt$ as a $W_t$-module; 
\item a non-trivial ideal contained in\/ $\ann_{W_t}\rsc FG$ and such
that the quotient\/ $W_t/\annFoG\simeq S'$~is holonomic.
\end{enumerate}
\end{cor}
\begin{proof}
{}From~\eqref{eq:3-formulations},
\begin{multline}\label{eqn:AnnWt}
\annFoGprime=\left(\annFxG\right)\cap\Wt\\
=\left(\FplusG\right)\cap\Wt,
\end{multline}
and we have shown (1)~in the corollary statement.  The $W_t$-module
$S'\simeq W_t/\annFoG$ is a holonomic $W_t$-module, as it is a
submodule of the holonomic $W_t$-module~$S$.  Now since $W_t$~is not
holonomic, $\annFoG$ must be non-trivial by a simple dimension
argument. Finally, we recall that this non-trivial ideal is contained
in $\ann_{W_t}\rsc FG$, since there is a surjection from~$S'$ to
$W_t/\ann_{W_t}\rsc FG$ given by $\psi:(u\otimes v)\mapsto \rsc
uv$. This proves (2)~in the corollary statement.
\end{proof}

\subsection{Termination}\label{sec:termination}

We now link the modules $S$ and~$S'$ to the algorithms and prove their
termination.  The termination of Algorithm~3 is more technical to
prove than that of Algorithm~1 since $\partial_{t_i}$~can act
separately on $F$ and~$G$.  Thus, for ease of presentation, we
consider Algorithms 1 and~3 in turn, to show that they eventually
generate a Gr\"obner basis for~$\annFoGt$.

\subsubsection{Termination of Algorithm~1}\label{sec:term_1}

The basic idea of Algorithm~1 is to compute filtrations of $\cI_F$
and~$\cI_G$ independently and incrementally and to recombine them at
each step.  The algorithm terminates when
condition~\eqref{item:breaking-condition-in-algo-1} in the algorithm
description is satisfied. We show that the algorithm will satisfy this
condition by eventually producing a Gr\"obner basis for
\mbox{$\annFoGt$}.  This subideal describes $F^\sadj\otimes G$
and~$\rsc FG$ as D-finite.

\begin{proof}(Theorem~\ref{thm:main}, Algorithm~1)
Algorithm~1 places a constraint on~$F$ that allows us to take
advantage of the simpler $W_t$-structure of~$U=W_{p,t}\cdot F$: since
each~$\partial_{t_i}\cdot F$ is~0, we have
$U=\field[t]\otimes_\field(W_p\cdot F)$
and~$I_F=\partial_tW_{p,t}+W_tJ_F$.  Taking the intersection
with~$\Wt$ is then far more transparent: from the previous section, we
obtain the following simplification of Eq.~\eqref{eqn:AnnWt}:
\begin{equation}\label{eqn:simplified}
\annFoG=\left(J_F^\sadj W_t+I_G\right)\cap W_t.
\end{equation}
Considering the monoid of monomials generated by
$p,\partial_p,\partial_t$, ordered by the monomial order~$\preceq$
specified by the algorithm, we denote by $\cV_\beta$ the filtration
$\bigoplus_{\gamma\preceq\beta}\field(t)\gamma$.

Assume that Algorithm~1 fails to terminate on some input $F$ and~$G$.
For any~$\beta$, Algorithm~1 thus eventually reaches a value for the
main loop index~$\alpha$ such that all the monomials that have been
considered in the algorithm span a vector space
containing~$\cV_\beta$.  After Step~\eqref{item:after-intro-in-algo-1}
in the main loop for this value~$\alpha$ of the loop index,
$B$~generates a vector space containing
\begin{equation*}
L_\beta:=
\bigl(J_F^\sadj W_t(t)\cap\cV_\beta\bigr)+\bigl(\cI_G\cap\cV_\beta\bigr).
\end{equation*}
By our choice of elimination term order, $B\cap W_t(t)$ consists of
generators of a vector space which contains the
intersection~$L_\beta\cap W_t(t)$.

Next, for each~$\gamma$, $\bigl(J_F^\sadj
W_t(t)+\cI_G\bigr)\cap\cV_\gamma$ is a subspace of~$L_\beta$ for
some~$\beta$.  Indeed, since $\cV_\gamma$ is finite-dimensional, so is
the intersection under consideration. Let us introduce a basis
$b_1,\dots,b_d$ of~it.  Each~$b_i$ can be written in the
form~$f_i+g_i$ for $f_i\in\cI_F^\sadj=J_F^\sadj W_t(t)$ and
$g_i\in\cI_G$, so that, provided $\beta = \max\{\max_i\deg f_i,
\max_i\deg g_i\}$, the intersection
\[
\bigl(J_F^\sadj W_t(t)+\cI_G \bigr)\cap\cV_\gamma
=\bigoplus_{i=1}^d\field(t)(f_i+g_i)
\]
is a subspace of
\[
\sum_{i=1}^d\field(t)f_i+\sum_{i=1}^d\field(t)g_i
\subset\bigl(W_t(t)J_F^\sadj\cap\cV_\beta\bigr)+\bigl(\cI_G\cap\cV_\beta\bigr)
=L_\beta.
\]

Since $\annFoGt$ is finitely generated by noetherianity of~$W_t(t)$,
we can choose a finite set of generators for it, and set~$\gamma$ to
their maximal leading monomial.  Consequently, the chosen generators
are in
\begin{equation*}
\annFoGt\cap\cV_\gamma=\bigl(W_t(t)J_F^\sadj+\cI_G\bigr)\cap
W_t(t)\cap\cV_\gamma.
\end{equation*}
By the reasoning above, the latter is a
subspace of~$L_\beta$ for some~$\beta$, and when the loop index
reaches a sufficiently high~$\alpha$, $\annFoGt$~is a subideal of the
ideal generated in~$W_t(t)$ by~$B\cap W_t(t)$.  Since, by
Corollary~\ref{thm:elements_of}, $W_t/\annFoG$ is a holonomic module,
$\annFoGt$~is of dimension~0, and
condition~\eqref{item:breaking-condition-in-algo-1} is satisfied.  The
algorithm terminates, a contradiction to our assumption.
\end{proof}

A limitation of the algorithm is that we cannot predict in advance how
many monomials must be tested, and hence cannot estimate the running
time.

\subsubsection{Termination of Algorithm~3}
\label{sec:term_3}
The termination of Algorithm~3 can be proved similarly, but we must
use greater care when treating the $\partial_{t_i}$.

\begin{proof}(Theorem~\ref{thm:main}, Algorithm~3)
Since there is no adjoint action for~$\partial_{t_i}$, we consider
occurrences of $\partial_{t_i}$ in the left argument of the scalar
product differently from those on the right side. This is modelled
in~$S$ by tensoring over $W_p[t]$, where $\partial_t$ is absent and
thus, $\partial_{t_i}\otimes 1$ differs from $1\otimes\partial_{t_i}$.
Both still obey the same commutation law with~$t_i$
as~$\partial_{t_i}$.  Denote the former by~$\partial_{\ell_i}$ and the
latter~by~$\partial_{r_i}$.

Having distinguished these two cases, we rewrite several of the
important elements from the previous proof using this new
notation. For example,
\begin{multline*}
\WoW=\field\bigl\langle
p,t,\partial_p,\partial_\ell,\partial_r;
[\partial_{p_i},p_j]=[\partial_{\ell_i},t_j]=[\partial_{r_i},t_j]=\delta_{i,j},\\
[p_i,p_j]=[p_i,t_j]=[t_i,t_j]=[\partial_{\ell_i},p_j]=[\partial_{r_i},p_j]=
[\partial_{p_i},t_j]=0
\bigr\rangle,
\end{multline*}
and its subalgebra~$\Wt$ is generated by $\field[t]$,
$\partial_{\ell_1}+\partial_{r_1}$,~\dots,
$\partial_{\ell_k}+\partial_{r_k}$.  We can also rewrite $\FplusG$ in the
form $\dlAnnF\field[\partial_r]+\field[\partial_\ell]\drAnnG$.
Algorithm~3 actually computes with coefficients that are rational
functions in~$t$, and so with elements of
$\dlAnnFt\field[\partial_r]+\field[\partial_\ell]\drAnnGt$.

In order to endow $\WoW$ with a filtration, let us extend the
ordering~$\preceq$ to monomials in
$p,\partial_p,\partial_\ell,\partial_r$ by considering any ordering
which, after setting $\partial_\ell=\partial_t,\partial_r=1$ or
$\partial_r=\partial_t,\partial_\ell=1$, respectively, induces the
ordering~$\preceq$.  We denote the extended ordering by~$\preceq$ as
well.  Then, we let $\cU_\beta$ denote the filtration
$\bigoplus_{\gamma\preceq\beta}\field(t)\beta$ for $\beta,\gamma$
ranging over the monomials in the variables
$p,\partial_p,\partial_r,\partial_\ell$.  Turning our attention
to~$\Wt(t)$, let $\cV'_\beta$ be the image of the~$\cV_\beta$ of the
previous section, under the same transformation which takes~$W_t(t)$
to~$\Wt(t)$, that~is,
\[
\cV'_\beta=
\bigoplus_{p^a\partial_p^b\partial_t^c\preceq\beta}
	\field(t)p^a\partial_p^b\left(\partial_\ell+\partial_r\right)^c.
\]
For each~$\beta$, there is~$\beta'$ such that
$\cV'_\beta\subset\cU_{\beta'}$.

Assume that Algorithm~3 fails to terminate on some input $F$ and~$G$.
Since the main loop enumerates all monomials in
$p,\partial_p,\partial_\ell,\partial_r$ in some order, for any~$\beta$
there exists a value of the index loop~$\alpha$ such that when the
loop reaches it, all monomials that have been enumerated span a vector
space containing~$\cU_\beta$.  After the algorithm has introduced
(variants of) $\alpha_F$ and~$\alpha_G$ at
Step~\eqref{item:after-intro-in-algo-2} for this value of~$\alpha$,
let us call $V_\alpha$ the vector space generated by the set~$B$.
Setting $\partial_\ell=\partial_t-\partial_r$ maps~$V_\alpha$ to a vector
space which contains
\[
H_\beta:=
\left(\dlAnnFt\field[\partial_r]\right)\cap\cU_\beta+
\left(\field[\partial_\ell]\drAnnGt\right)\cap\cU_\beta.
\] 
We use this fact to conclude termination.

At this point we show that for each~$\gamma$, the vector space
$\cX\cap\cV'_\gamma$ where
\[\cX=\FplusGt\]
is a subspace of~$H_\beta$ for some~$\beta$.  Indeed, choose~$\gamma'$
such that $\cV'_\gamma\subset\cU_{\gamma'}$, so that
$\cX\cap\cV'_\gamma\subset\cX\cap\cU_{\gamma'}$.  The latter
intersection is finite-dimensional, since $\cU_{\gamma'}$~is~so.
Suppose it has for basis $b_1,\dots,b_d$, with each~$b_i$ of the form
$b_i=f_ir_i+l_ig_i$, where $f_i\in\dlAnnFt$, $g_i\in\drAnnGt$,
$r_i\in\field[\partial_r]$, and $l_i\in\field[\partial_\ell]$, and set
$\beta=\max\{\max_i\deg f_ir_i,\max_i\deg l_ig_i\}$, where here
$\deg$~extracts the leading monomial.  Then,
\[ 
\cX\cap\cV'_\beta\subset
\bigoplus_{i=1}^d\field(t)(f_ir_i+l_ig_i)\subset
\sum_{i=1}^d\field(t)f_ir_i+\sum_{i=1}^d\field(t)l_ig_i\subset
H_\beta.
\]

By noetherianity, we can choose a finite set of generators for
$\annFoGt$, and set~$\gamma$ to their maximal leading monomial.  The
generators are thus elements of $\annFoGt\cap\cV_\gamma$, which is
isomorphic to $\annFoGtprime\cap\cV'_\gamma$.  By~\eqref{eqn:AnnWt}
the latter is also $\cX\cap\cV'_\gamma$, and, as explained above,
there is~$\beta$ such that this is a subspace of~$H_\beta$.

By our earlier loop invariant, the same generators, after setting
$\partial_\ell=\partial_t-\partial_r$, are contained in the space spanned
by~$B$ when the loop index reaches a sufficiently high~$\alpha'$.
Thus, it suffices to run the algorithm until this~$\alpha$ and
generators of $\annFoG$ will be contained in~$B$.  At this point the
termination conditions are satisfied, and the algorithm terminates.
\end{proof}

\section{Asymptotic Estimates}\label{sec:asympt}

We now illustrate how the differential equations computed by our
algorithms may be exploited in order to derive asymptotic estimates
of combinatorial quantities.

\subsection{Outline of the method}
A very general principle in asymptotic analysis is that the asymptotic
behaviour of a sequence is governed by the local behavior of its
generating series at its singularity of smallest modulus, see for
instance~\cite[Section~10]{Odlyzko95}.  Our approach is thus based on
applying the classical analysis of linear differential equations as
presented in textbooks such as~\cite{Ince56,Wasow87} in order to
derive asymptotic estimates for the coefficients. Moreover, large
parts of this analysis can be automated thanks to the algorithms
described in~\cite{Malgrange79,Tournier87,vanHoeij97c}, many of which
have been implemented in computer algebra systems\footnote{In Maple,
this functionality is provided by {\tt DEtools[formal\_sol]}.}. An
alternative approach based on Birkhoff's work can be found
in~\cite{WiZe85}.

In the special case of solutions of linear differential equations, the
possible location of singularities is restricted to the roots of the
coefficient of the highest derivative. Then, the analysis depends on
the nature of the singularity. The classical theory distinguishes two
kinds of singular points: regular singular points, where the solutions
have an algebraic-logarithmic behavior; and irregular singular points
where the solutions have an essential singularity of the type
exponential of a rational power. Accordingly, the asymptotic behavior
of the coefficients is deduced either by singularity
analysis~\cite{FlOd90b,Jungen31}, or by the saddle-point
method~\cite{Hayman56,Wyman59}; both approaches are implemented in the
{\tt algolib} library.

This asymptotic analysis of D-finite generating series extends to the
divergent case. Indeed, the coefficients~$u_n$ of a divergent D-finite
series grow at most like a power of~$n!$ with a rational
exponent~$p/q$ which can be computed (see example below). Then one
constructs an auxiliary differential equation satisfied by the
convergent generating series of~$u_n/(n(n-q)(n-2q)\dotsm r)^p$ (where
$r$~denotes the remainder of the division of~$n$ by~$q$), to which the
previous method applies. This construction is achieved thanks to the
closure properties of D-finite series, by multiplying~$u_n$ with the
solution of the recurrence~$(n+q)^pv_{n+q}=v_n$, which, up to a
constant, grows like~$n!^{p/q}n^{p(q-1)/2q}$. This operation is
implemented in the {\tt gfun} package.

\subsection{$k$-uniform Young tableaux}
We now illustrate this method in the special case of the $k$-uniform
Young tableaux of Section~\ref{sec:young}. We treat in detail the
case~$k=3$; other cases are similar. To the best of our knowledge,
these asymptotic estimates are new.

We start from the differential equation for~$k=3$ to be found in
Table~\ref{ftab}. This is a second-order differential equation and its
leading coefficient vanishes at the origin. This indicates a possible
singularity of~$Y_3(t)$ at the origin, which would be reflected by the
divergence of this series. Indeed, from this differential equation, a
linear recurrence is readily computed for the
coefficients~$u_n:=y_n^{[3]}$:
\begin{multline*}
u_n+u_{n+1}-(3n+12)u_{n+2}-4u_{n+3}+(6n+35)u_{n+4}-15u_{n+5}\\
+(9n^2+93n+242)u_{n+6}+(18n+126)u_{n+7}-(9n^2+159n+698)u_{n+8}\\
+(9n^2+147n+606)u_{n+9}-(18n^2+366n+1884)u_{n+10}\\
-(48n+552)u_{n+11}+(24n+288)u_{n+12}=0.
\end{multline*}

\subsubsection{Divergence}
{}From this recurrence it is easy to compute a couple hundred
coefficients and observe their rapid growth. Simple experiments
indicate that the growth of these coefficients is of
order~$\sqrt{n!}$. That this growth is the exact exponent of~$n!$ in the
behavior follows upon considering the degrees of the coefficients in
the recurrence: the terms of order~12 and~11 have coefficients of
degree~1, while the term of order~10 has a coefficient of degree~2
(the maximal degree). Thus, up to first order, the behavior is dictated by
\[24nu_{n+12}=18n^2u_{n+10},\]
which leads to a growth of order~$(\frac34)^{n/2}n!^{1/2}$. 
In order to derive a more precise estimate, we compute a linear
differential equation satisfied by the \emph{convergent} generating
function 
of~$y_n^{[3]}v_n$ where $v_n$
satisfies~$v_{n+2}=v_n/(n+2)$. This differential 
equation is obtained by first computing a linear recurrence
for~$y_n^{[3]}v_n$, which exists thanks to the closure properties
of linear recurrent sequences. This closure operation produces
a linear recurrence of order~24 with coefficients of degree~29. From
there we obtain a linear differential equation of order~29 with
coefficients of degree~37, which we now analyze. 

\subsubsection{Singular behavior} The leading coefficient of the
previous equation is~$t^{27}(3t^{2}-4)$, up to a constant factor. This
reveals a dominant singularity at~$\rho=2/\sqrt{3}$, thus confirming
the growth order~$(3/4)^{n/2}$ expected from the previous
stage\footnote{We could also have incorporated this factor in the
recurrence for~$v_n$.}. The next step consists in analyzing the
behavior of our convergent generating series in the neighborhood
of~$\rho$. A local analysis of the
differential equation reveals that all solutions of
this equation of order~29 behave like 
\[g(u)+\lambda \frac{\exp\left(\frac{3}{4u}\right)}{\sqrt{u}} \left(1-{\frac {145}{144}}u-{\frac {8591}{41472}}{u}^{2}+O \left( {u}^{3} \right)\right),\qquad 1-z/\rho=u\rightarrow0,\]
where~$g$ is an analytic function at~0, and $\lambda$ is a constant
depending on the solution.

\subsubsection{Asymptotic estimate} This behavior is typical of an
irregular singular point and can thus be dealt with using the
saddle-point method. Putting everything
together, we finally obtain 
\[y_n^{[3]}=C_3n!^{1/2}\left(\frac{\sqrt{3}}{2}\right)^n\frac{\exp{\sqrt{3n}}}{n^{3/4}}(1+O(1/n)),\]
for some constant~$C_3$, and where the $O$-term hides the beginning of an expansion in descending powers of~$n$ that could be computed with the same method.

The constant~$C_3$ can then be approximated numerically by using
Romberg's acceleration method, adapted to powers of~$n^{-1/2}$, and we
get:
\[C_3\approx0.377200.\]

\subsubsection{Other values of~$k$}
\begin{table}
\begin{center}
\begin{tabular}{lcc}\hline\\
1&$\displaystyle C_1
\frac{\exp{\sqrt{n}}}{\sqrt{n!}\,n^{1/4}}$&$C_1\approx0.347829$\\[3mm]
2&$\displaystyle C_2 \frac{\exp{\sqrt{2n}}}{\sqrt{n}}$&$C_2\approx0.282094$\\[3mm]
3&$\displaystyle C_3 \sqrt{n!}\left(\frac{\sqrt{3}}{2}\right)^n\frac{\exp{\sqrt{3n}}}{n^{3/4}}$&$C_3\approx0.377200$\\[4mm]
4&$\displaystyle C_4 n!\left(\frac{2}{3}\right)^n\frac{\exp{2\sqrt{n}}}{n}$&$C_4\approx0.831565$\\[4mm]
\hline\\
\end{tabular}
\caption{Asymptotic number of~$k$-uniform Young tableaux}
\label{table-asympt}
\end{center}
\end{table}

The computation of the asymptotic behavior of~$y^{[k]}_n$ for other
values of~$k$ is completely similar, provided one has computed the
differential equation. We summarize our results in
Table~\ref{table-asympt}.  This serves to illustrate a typical use of
our techniques in experimental mathematics to obtain conjectures such
as the following.
\begin{conj}The number~$y_n^{[k]}$ of $k$-uniform Young tableaux
of size~$n$ behaves asymptotically according to
\begin{equation*}
y_n^{[k]}\sim
\frac1{\sqrt2}\left(\frac{e^{k-2}}{2\pi}\right)^{k/4}
n!^{k/2-1}\left(\frac{k^{k/2}}{k!}\right)^n\frac{\exp(\sqrt{kn})}{n^{k/4}},
\qquad n\rightarrow\infty.
\end{equation*}
\end{conj}

This conjecture is proved for $k=1$ and~$k=2$: the constant is
obtained from a closed form solution of the differential equation.
For $k=3$ and~$k=4$, only the value of the constant is conjectural.
The proof of the general case of the conjecture requires techniques
such as those of~\cite{GoMc90,McKay90}, which fall outside of the
scope of this article.

\subsection{Conclusion}
The main advantages of our method are its general applicability, its
ability to produce full asymptotic expansions up to \emph{one\/}
constant factor, the availability of computer algebra programs that
automate many of its steps. The price to pay for this generality is
that the method can only produce numerical estimates for the constant
factor. In some special cases, specific approaches often exist that
provide this constant term.

\section{Conclusions and Directions for Future Work}\label{sec:concl}

\subsection{Applying the method to other scalar products}
\label{sec:other-scalar-products}

Let us note that the method of this article can be applied in the case
of other scalar products, provided that the corresponding
adjunction~$\sadj$ (no longer denoting the symmetric adjunction) is a
linear involution that preserves the total degree (in~$p,\partial_p$)
of the differential operators.  In effect, one should simply set
$M=(U^\sadj)^\uadj$ and~$N=V$ to obtain a suitable analogue
to~\eqref{eq:adj-fact} and prove the holonomy, thus D-finiteness, of
the scalar product: $M$~is holonomic if and only if $U$~is.  Since the
statement and proof of Algorithm 1 and~3 do not make use of any other
special property of~$\sadj$ than being a degree-preserving involution,
correctness of the algorithms can then be established along the same
lines as for the case of the scalar product of symmetric 
functions. 

We use this idea in the next two sections by introducing various
scalar products given by an adjunction relation involving a formal
parameter.

\subsection{Calculating the Kronecker product of symmetric functions}
\label{sec:kronecker-product}

Another symmetric function operation, closely related to the scalar
product, is the Kronecker product, also known as the tensor
product. One can define it on the power basis as $p_\lambda*p_\mu
=\rsc{p_\lambda}{p_\mu}p_\lambda$. Gessel showed in~\cite{Gessel90}
that given two D-finite symmetric series $F$ and~$G$, the Kronecker
product $F*G$ is also a D-finite symmetric series. Algorithm~1 can be
used to make this fact effective via the following observation:
\[
p_\lambda*p_\mu =\rsc{p_\lambda t^\lambda}{p_\mu}\big|_{t_i=p_i}.
\]
More precisely, we rewrite a Kronecker product as a scalar product by
multiplying each~$p_i$ in~$F$ by~$t_i$. In the system which results we
make the substitution $t_i=p_i$ and $\partial_{t_i}=\partial_{p_i}$.

We formalize this in the following algorithm, which merely calls
Algorithm~1 on modified input systems.
\begin{algo}[Kronecker Product]\label{thm:algo4}
\mbox{}\\
{\sc Input:} Symmetric functions\/ $F\in\field[[p]]$ and\/
$G\in\field[[p]]$, both D-finite in\/~$p$, each given by a D-finite
description in\/~$W_p$.\\
{\sc Output:} A D-finite description of\/~$F*G$ in\/~$W_t$.
\begin{enumerate}

\item Call~$\cG$ the system defining~$G$ and set
$\cG'=\{t_1\partial_{t_1}-p_1\partial_{p_1},\dots,t_n\partial_{t_n}-p_n\partial_{p_n}\}$;

  \begin{enumerate}

  \item For each element in~$\cG$, replace $p_i$ with $t_ip_i$,
  $\partial_{p_i}$ with $t_i^{-1}\partial_{p_i}$ and add to~$\cG'$;

  \item For each element in~$\cG$, replace $p_i$ with $t_ip_i$,
  $\partial_{p_i}$ with $p_i^{-1}\partial_{t_i}$, clear denominators,
  and add to~$\cG'$;

\end{enumerate}

\item Follow the steps of Algorithm~1 on the input system for~$F$ and
the modified system~$\cG'$ for~$G$;

\item In the output of Algorithm~1 make the substitution $t_i=p_i$ and
$\partial_{t_i}=\partial_{p_i}$ and return this value.

\end{enumerate}
\end{algo}
Many interesting problems which use this operation require an infinite
number of $p_n$, and are thus at first glance seemingly unsuitable for
direct application of our algorithms. However, applying our algorithms
for several truncations of a combinatorial problem can serve as a
means to generate information upon which reasonable conjectures can be
formulated. For example, Eq.~\eqref{eqn:schur} below was initially
conjectured after a clear pattern emerged from a sequence of appeals
to Algorithm~\ref{thm:algo4}. For each of these, we render the problem
applicable by setting most $p_n$'s to 0. In some cases, notably
symmetric series arising from plethysms, there is sufficient symmetry
and structure which can be exploited to verify these guesses by
applying one of Algorithm~4 to well chosen subproblems.  That is, in
certain cases, such as the example that follows, the Kronecker product
of two functions each with an infinite number of $p_n$ variables can
be reduced to a finite number of symbolic calculations.

For example, if two symmetric series $F$ and $G$ can be expressed
respectively  in the form 
\[F(p_1, p_2, \ldots)=\prod_{n\ge 1} f_n(p_n) \qquad \text{ and } \qquad
G(p_1, p_2, \ldots)=\prod_{n\ge 1} g_n(p_n),\]
for functions $f_n$, $g_n$, then one can
easily deduce that 
\begin{equation}\label{eqn:prodform}
F*G=\prod_{n\ge 1}  f_n(p_n)*g_n(p_n).
\end{equation}
Remark that series which arise as plethyms of the form $h[u]$ or
$e[u]$, where $u$ can be written as a sum $\sum_n u_n(p_n)$, for some
functions $u_n$, are precisely of this form. For example, we can use
this fact to compute the Kronecker product of the sum of all Schur functions
\[F(p_1, p_2, \ldots)=\sum_\lambda
s_\lambda=h[p_1+1/2p_1^{2}-1/2p_2]=\exp\left(\sum_i\frac{p_i^2}{2i}+\frac{p_{2i-1}}{2i-1}\right),\]
and itself.  Due to the patterns present, we can
reduce the calculation of the entire product to two symbolic
calculations. More precisely, in order to determine a system of
differential equations satisfied by $G=F*F$ we consider only the even
and odd cases, and set 
\[f_{2n}=\exp(p_{2n}^2/{4n})\quad \text{ and }\quad
f_{2n-1}=\exp((p_{2n-1}^2/2+p_{2n-1})/(2n-1)).
\]
All of the functions $g_{2n}=f_{2n}*f_{2n}$ are obtained from a single
computation by our Algorithm~4, adapted to handle a formal
parameter. This modification is of the same nature of that described
in Section~9.1. Here we introduce the scalar product given by the
adjunction formula $p^\diamond=n\partial$ for a {\em formal
parameter\/} $n$ from the field~$K$. Thus computing
$\exp(p^2/4n)*\exp(p^2/4n)$ with this variant algorithm results in a
first-order operator in $p$ and $\partial$, which, once interpreted
back in terms of $p_n$ becomes:
\begin{equation*}
(1-p_n^2)\frac{\partial g_n(p_n)}{\partial p_n}+p_ng_n(p_n)=0,
\qquad\text{for even~$n$}.
\end{equation*}
A second calculation for $g_{2n-1}=f_{2n-1}*f_{2n-1}$ results in:
\begin{equation*}
n(1+p_n)(1-p_n)^2\frac{\partial g_n(p_n)}{\partial p_n}-
\left(1+(n+1)p_n-np_n^2\right)g_n(p_n)=0,
\qquad\text{for odd~$n$}.
\end{equation*}
These linear equations are satisfied respectively by the functions
\begin{multline*}
g_{2n}=\left(1-p_{2n}^2\right)^{-1/2}
\qquad\text{and}\\
g_{2n-1}=\exp\left(\frac{p_{2n-1}}{(2n-1)(1-p_{2n-1})}\right)\left(1-p_{2n-1}^2\right)^{-1/2}.
\end{multline*}
Applying Eq.~\eqref{eqn:prodform} above, we get the following result.

\begin{prop}
The Kronecker product of the sum of the Schur functions with itself is
\begin{equation}\label{eqn:schur}
\left(\sum_\lambda s_\lambda\right)*\left(\sum_\lambda
s_\lambda\right)={\exp\left(\sum_{n\geq
1}\frac{p_{2n-1}}{(2n-1)(1-p_{2n-1})}\right)}{\left({\prod_{n\geq
1}\left(1-p_n^2\right)}\right)^{-1/2}}.
\end{equation}
\end{prop}

\subsection{A $q$-analogue}

A $q$-calculus parameter can be incorporated in symmetric functions in
several~ways.

Apart from the scalar product defined by~\eqref{eq:scalhm}, several
other ones are of interest in relation to symmetric functions, notably
the following two, which lead to the definitions of Hall and Macdonald
polynomials respectively:
\[
\rsc{p_\mu}{p_\lambda}=z_\lambda\delta_{\mu,\lambda}
    \prod_{i=1}^{l(\lambda)}(1-t^{\lambda_i})
\quad\text{ and }\quad
\rsc{p_\mu}{p_\lambda}=z_\lambda\delta_{\mu,\lambda}
    \prod_{i=1}^{l(\lambda)}\frac{(1-t^{\lambda_i})}{1-q^{\lambda_i}},
\]
where $\ell(\lambda)$~is the length~$k$ of a partition
$\lambda=(\lambda_1,\dots,\lambda_k)$.  The same approach as in this
article works in this setting and our Maple code has been adapted very
easily\footnote{This variant is also available at {\tt
http://algo.inria.fr/mishna}.}.

As a related problem, the ring homomorphism
$\theta_q:\Lambda\rightarrow\field[q][[t]]$ defined as
\begin{equation*}
\theta_q\bigl(f(x_1,x_2,\ldots)\bigr)
=f\bigl((1-q)t,(1-q)qt,(1-q)q^2t,\ldots\bigr)
\end{equation*}
is useful for studying partitions and for counting permutations
\cite{Stanley99}.  This is one possibility for a $q$-analogue to the
map~$\theta$ from Theorem~\ref{thm:theta} (named exponential
specialization in~\cite{Stanley99}), since $\lim_{q\rightarrow
1}\theta_q(F)=\theta(F)(x)$.  An algorithm to compute~$\theta_q$,
possibly mapping differential equation to $D_q$ equation should be of
interest.

\subsection{Other conditions for D-finite closure}

Remark that Theorem~\ref{thm:rsc_pres_df} requires that $g$~be a
function of only a finite number of~$p_n$.  The necessity of this
condition is evident in the following example.  Find a sequence~$c_n$
such that $\sum c_n t^n$ is not D-finite.  However, according to the
given definition of D-finite symmetric series, $\sum_n c_n p_n$ is
D-finite, as is $\sum_n p_nt^n/n$.  The series $\rsc{\sum_n c_n
p_n}{\sum_n p_nt^n/n}=\sum_n c_n t^n$ is not D-finite by construction.

On the other hand, the condition is not essential.  We have that
$\rsc{H(1)}{H(t)}=\frac{1}{1-t}$, which is D-finite despite $H$~being
a function of {\em all $p_n$.}  Perhaps a closer investigation on the
level of modules could reveal a refined condition.

\bigskip
\subsection*{Acknowlegements}
The authors wish to thank Fran\c cois Bergeron for promoting D-finite
symmetric functions as an interesting area of study. The second author
also extends gratitude towards NSERC for funding, and to Projet {\sc
Algo}, {\sc Inria}, for their generous invitations during which much
of the work was completed.  Finally, we thank the anonymous referees
who read the work carefully and offered many useful suggestion to
improve the clarity.

\appendix

\section{4-Uniform Young Tableaux}\label{sec:fourreg}

The differential equation satisfied by $Y_4(t)$ is
\begin{multline*}
64t^4(t-2)^2(t+1)^4\alpha(t)Y_4^{(3)}(t)
-16t^2(t-2)(t+1)^2\beta(t)Y_4^{(2)}(t)\\
+4\gamma(t) Y_4'(t)
-\delta(t) Y_4(t)=0
\end{multline*}
where~$\alpha(t),\beta(t),\gamma(t),\delta(t)$ are irreducible
polynomials given by
\begin{align*}
\alpha(t)&=t^{14}-t^{13}-5t^{12}-7t^{11}+6t^{10}
+35t^9+39t^7-50t^6-162t^5-92t^4\\
&+228t^3+424t^2+248t+48,\\
\beta(t)&=t^{29}-3t^{28}-16t^{27}+24t^{26}+147t^{25}
+14t^{24}-770t^{23}-666t^{22}+1416t^{21}\\
&+3567t^{20}-916t^{19}-16598t^{18}+17766t^{17}+40678t^{16}
-102556t^{15}\\
&-53272t^{14}+390656t^{13}+364080t^{12}-707936t^{11}-1406336t^{10}
-552544t^9\\
&+1397664t^8+2020864t^7+176256t^6-916864t^5+304896t^4+1283328t^3\\
&+877056t^2+253440t+27648,\\
\gamma(t)&=t^{28}-t^{27}-14t^{26}-20t^{25}+111t^{24}
+278t^{23}-196t^{22}-1216t^{21}\\
&-1384t^{20}+2765t^{19}+3170t^{18}-3400t^{17}+12140t^{16}+15588t^{15}\\
&-70280t^{14}-108946t^{13}+121796t^{12}+349056t^{11}+116992t^{10}-481704t^9\\
&-706320t^8+3040t^7+581184t^6+158688t^5-297408t^4-173952t^3\\
&+22272t^2+35712t+6912,\\
\delta(t)&=2t^{21}-3t^{20}-17t^{19}-2t^{18}+74t^{17}+105t^{16}-108t^{15}-172t^{14}-252t^{13}\\
&+432t^{12}-667t^{11}+1500t^{10}+7336t^9-3772t^8-23056t^7-20584t^6\\
&+15504t^5+38160t^4+17904t^3-4512t^2-5568t-1152.
\end{align*}

\section{Sample Maple Session for 3-Regular Graph Computation}

The following Maple session indicates the user-level routines required
to program Algorithm~2. It requires the library {\tt algolib}, which
is available at
\verb+http://+\discretionary{}{}{}\verb+algo.inria.fr/packages/+.

\begin{verbatim}
# Load the packages.
with(Ore_algebra): with(Mgfun): with (Groebner):
# Determine the DE satisfied by the generating function
# for 3-regular graphs.
k:=3: Fp:= exp(1/2*p1^2-1/4*p2^2-1/2*p2+p3^2/6):
Gp:=exp(1/6*t3*p1^3+1/2*t2*p1^2+t1*p1+1/2*t3*p2*p1
    +1/2*t2*p2+1/3*t3*p3):
# Define the variables.
vars:= seq(p||i, i=1..k):  dvars:= seq(d||i, i=1..k):
tvars:= seq(t||i, i=1..k): dtvars:= seq(dt||i, i=1..k):

# Define the algebra.
A:= diff_algebra(seq([dvars[i], vars[i]], i=1..k), 
seq([dtvars[i], tvars[i]], i=1..k), polynom={vars}):
At:= diff_algebra(seq([dtvars[i], tvars[i]], i=1..k)):

# Define the monomial orders.
T[g]:=termorder(A, lexdeg([dvars, vars],[dtvars])): 
T[f]:=termorder(A,tdeg(vars, dvars, dtvars)):

# Define the systems.
sys[g]:=dfinite_expr_to_sys(Gp, F(seq(p||i::diff, i=1..k),
        seq(t||i::diff, i=1..k))):
newsys[g]:=subs(
    [seq(diff(F(vars,tvars),vars[i])=dvars[i],i=1..k),
     seq(diff(F(vars, tvars), tvars[i])=dtvars[i], i=1..k), 
     F(vars,tvars)=1], sys[g]):

# Find the Groebner basis for G.
GB[g]:=gbasis(newsys[g],T[g]);

# Do the same for F.
sys[f]:=dfinite_expr_to_sys(Fp, F(seq(p||i::diff, i=1..k))):
newsys[f]:=subs([seq(diff(F(vars),vars[i])=dvars[i],i=1..k),
F(vars)=1],sys[f]); 
GB[f]:=gbasis(newsys[f],T[f]);

# Define the adjoint and reduction procedures.
star:= x->subs(
    [seq(d||i=1/i*p||i, i=1..k),seq(p||i=d||i*i, i=1..k)],x):
rdc[f] := x->star(star(x)-map(normalf, star(x), GB[f], T[f]));
rdc[g] := x->normalf(x, GB[g], T[g]);

# Reduce the Groebner basis of F.
for pol in GB[f] do m[pol]:=rdc[g](pol) end do:

# Small optimization: we will always try to reduce with respect
# to a linear term when possible.
lpol:=[seq(m[i],i=subsop(1=NULL,GB[f])),m[GB[f][1]]]:

for indelim from k-1 by -1 to 1 do 
    # eliminate dt.indelim
    for j from 2 to nops(lpol) do 
       newpol[j]:=skew_elim(lpol[j],lpol[1],dt||indelim,At)
    end do;
    # set t.indelim = 0
    lpol:=map(primpart,subs(t||indelim=0,
        [seq(newpol[j],j=2..nops(lpol))]),[dtvars])
end do:

# The only term left is the correct one.
ode:=op(lpol):
# Convert to recurrence.
REC:=diffeqtorec(
    {applyopr(ode, F(t||k), At), F(0)=1}, F(t||k), a(n)):
# Calculate some terms.
GRAPH:=rectoproc(REC, a(n),list)(20):
[seq(GRAPH(10)[i]*(i-1)!,i=1..20)];
 
 
 [1,0,0,0,1,0,70,0,19355,0,11180820,0,11555272575,0,
   19506631814670,0,50262958713792825,0,187747837889699887800,0]  
\end{verbatim}

\bibliographystyle{acm}
\bibliography{effective}

\end{document}